\documentclass{amsart}

\usepackage{amsmath,amsthm}
\usepackage{amsfonts,amssymb}

\newtheorem{thm}{Theorem}[section]
\newtheorem{cor}[thm]{Corollary}
\newtheorem{lem}[thm]{Lemma}
\newtheorem{prop}[thm]{Proposition}

\newtheorem{algo}[thm]{Algorithm}
\newtheorem{defn}[thm]{Definition}
 
\theoremstyle{remark}
\newtheorem{rem}{Remark}[section]

 \def\xb{{\mathbf x}}

 \def\CA{{\mathcal A}}     
 \def\CB{{\mathcal B}}
 \def\CD{{\mathcal D}}

 \def\CI{{\mathcal I}}
 \def\CJ{{\mathcal J}}
 
 \def\CO{{\mathcal O}}
      
 \def\CR{{\mathcal R}}     
      
 \def\CV{{\mathcal V}}

 \def\NN{{\mathbb N}}
 \def\PP{{\mathbb P}}
 
 \def\RR{{\mathbb R}}

        \def\proj{\operatorname{proj}}
        \def\sign{\operatorname{sign}}
        \def\supp{\operatorname{supp}}

\newcommand{\wt}{\widetilde}
\newcommand{\wh}{\widehat}

\begin{document}
\title[reconstruction of images from Radon projections]
{A new approach to the reconstruction of images from Radon projections} 
 
\author{Yuan Xu}
\address{Department of Mathematics\\ University of Oregon\\
    Eugene, Oregon 97403-1222.}
\email{yuan@math.uoregon.edu}

\date{\today}
\thanks{The author was supported in part by the NSF Grant DMS-0201669}

\begin{abstract}
A new approach is proposed for reconstruction of images from Radon 
projections. Based on Fourier expansions in orthogonal polynomials of two 
and three variables, instead of Fourier transforms, the approach provides 
a new algorithm for the computed tomography. The convergence of the algorithm 
is established under mild assumptions.
\end{abstract}

\maketitle

\section{Introduction}
\setcounter{equation}{0}

The fundamental problem of the computed tomography (CT) is to reconstruct
an image from its Radon projections  (X-rays).  To be more precise, let us 
consider the case of a two dimensional image, described by a density function 
$f(x,y)$ defined on the unit disk $B^2=\{(x,y):x^2+y^2 \le 1\}$ of the $\RR^2$ plane. 
A Radon projection of $f$ is a line integral,
$$
\CR_\theta(f;t) :=  \int_{I(\theta,t)} f(x,y) dx dy, 
\qquad 0 \le \theta \le 2\pi, \quad -1\le t \le 1, 
$$
where $I(\theta,t) = \{(x,y): x \cos \theta + y \sin \theta = t\} \cap B^2$ 
is a line segment inside $B^2$.  The reconstruction problem of CT requires 
finding $f$ from its Radon projections. This question is completely solved
if the complete data is given, that is, if projections for all $t$ and 
$\theta$ are given. In practice, however, only a finite number of projections 
can be measured. Hence, the essential problem of CT is to find an effective 
algorithm that produces a good approximation to $f$ based on a finite 
number of Radon projections.  
The arrangement of the projections 
is referred as scanning geometry  since it is determined by the design of 
the scanner.  

Currently, the most important method in CT is the filtered backprojection 
(FBP) method which, like several other methods, is based on techniques of
Fourier transforms.  In principle, the FBP method works 
for banded limited functions.  It requires choosing a function $W_b$ 
(a low-pass filter with cut-off frequency $b$) that approximates the 
$\delta$-distribution, and it includes steps of linear interpolation and a 
discrete convolution that approximates a continuous convolution.  Each
of these steps could introduce serious errors in the algorithm.  See the 
discussion in \cite{N, NW}.

The purpose of the present paper is to develop a direct approach
for reconstruction of images in CT. Instead of using Fourier 
transforms, we will use orthogonal expansions based on 
orthogonal polynomials on $B^2$. Let $S_n f$ denote the $n$-th 
partial sum of such an expansion. It turns out that $S_{2m} f$ satisfies 
the following remarkable formula
\begin{equation}  \label{eq:1.1}
S_{2m}(f;x,y) =  \sum_{\nu=0}^{2m} \frac{1}{\pi} 
         \int_{-1}^1 \CR_{\phi_\nu}(f;t) \Phi_\nu (t; x,y) dt,
\end{equation}
where  $\phi_\nu = 2 \nu \pi/(2m+1)$ are equally spaced angles 
along the circumference of the disk and $\Phi_\nu$ are polynomials 
of two variables given by explicit formulas (see Theorem 
\ref{thm:partial-sum}).  Applying appropriate quadrature formulas 
to the integrals in \eqref{eq:1.1} leads to an approximation, 
$\CA_{2m} f$, of $f$ that uses discrete Radon data.  For example, 
using an appropriate  Gaussian quadrature formula leads to an $\CA_{2m} f$ 
that uses the Radon data
\begin{equation}  \label{eq:1.2}
  \{\CR_{\phi_\nu}( f; \cos \tfrac{j\pi}{2m+1}) : \quad 0\le \nu \le 2m, \quad 
       1 \le j \le 2m \},
\end{equation}
of the scanning geometry of parallel beams, and the polynomial 
$\CA_{2m} f$ is given by an explicit formula 
\begin{equation} \label{eq:1.3} 
 \CA_{2m}(f; x,y) = \sum_{\nu=0}^{2m} \sum_{j=1}^{2m}
   \CR_{\phi_\nu}(f;\cos \tfrac{j\pi}{2m+1}) T_{j,\nu}(x,y),
\end{equation} 
where $T_{j,\nu}(x,y)$ are polynomials that are given by simple
formulas. The operator $\CA_{2m}f$ provides a direct 
approximation to $f$, which is the essence of our reconstruction
algorithm. The new algorithm can be implemented easily as the
$T_{j,\nu}$  are fixed polynomials that can be stored in a
table beforehand.  Furthermore, the construction 
allows us to add a multiplier factor and there is a natural 
extension of the algorithm to a cylindrical domain in $\RR^3$.  

The operator $\CA_{2m}$ in \eqref{eq:1.3} reproduces polynomials
of degree $2m -1$. Furthermore, we will prove that the operator
norm of $\CA_{2m}$ in the uniform norm is $\CO(m \log (m+1))$,
only slightly worse than the norm of $S_{2m}$, which is $\CO(m)$.
As a consequence, it follows that the algorithm converges uniformly
on $B^2$ if $f$ is a $C^2$ function.  There is no need to assume
that $f$ is banded limited.  

At the moment the author is working with Dr. Christoph Hoeschen
and Dr. Oleg Tischenko, Medical Physics Group of the National 
Research Center for Environment and Health, Germany, to implement
the algorithm numerically. The results that we have obtained so 
far are very promising and will be reported elsewhere.

The paper is organized as follows. The background of Radon
transforms and orthogonal polynomials is discussed in the next
section, where an identity that will play a fundamental
role in the analysis is also proved.  The  Fourier orthogonal 
expansions are studied in the
Section 3, where  the identity \eqref{eq:1.1} will be proved. 
The reconstruction algorithms are presented in Section 4.
The convergence of the basic algorithm for the 2-D reconstruction
is proved in Section 5 . 
 
\section{Radon transform and polynomials}
\setcounter{equation}{0}

Let $B^2= \{(x,y): x^2+y^2 \le 1\}$ denote the unit disk on the plane.
It is often convenient to use the polar coordinates 
$$
x= r \cos \theta, \quad y = r \sin \theta, \qquad r \ge 0, \quad
   0 \le \theta \le 2 \pi.
$$ 
Let $\theta$ be an angle as in the polar coordinates; that is, $\theta$ 
is measured counterclockwise from the positive $x$-axis. Let $\ell$ denote 
the line $\ell(\theta,t) = \{(x,y): x \cos \theta + y \sin \theta = t\}$
for $-1 \le t \le 1$. Clearly the line is perpendicular to the direction
$(\cos \theta,\sin\theta)$ and $|t|$ is the distance between the
line and the origin. We will use the notation
\begin{equation} \label{eq:line}
I(\theta, t) = \ell (\theta, t) \cap B^2, \qquad 
      0 \le \theta< 2 \pi, \quad -1 \le t \le 1,
\end{equation}
to denote the line segment of $\ell$ inside $B^2$. The points on $I(\theta,t)$
can be represented as follows: 
$$
x = t \cos \theta - s \sin \theta, \quad  y = t \sin \theta + s \cos \theta,
$$
for $s \in \left[-\sqrt{1-t^2}, \sqrt{1-t^2}\right]$.

The Radon projection of a function $f$ in the direction $\theta$ with
parameter $t \in [-1,1]$ is denoted by $\CR_\theta (f;t)$,   
\begin{align} \label{eq:Radon}
\CR_\theta(f;t) := & \int_{I(\theta,t)} f(x,y) dx dy \\
  = & \int_{-\sqrt{1-t^2}}^{\sqrt{1-t^2}} 
  f(t \cos\theta - s \sin \theta, t\sin \theta + s \cos \theta) ds. \notag
\end{align}
In the literature a Radon projection of a bivariate function is also called 
an $X$-ray. The definition shows that $\CR_\theta(f;t)=\CR_{\pi+\theta}(f;-t)$.

Let $\Pi^2$ denote the space of polynomials of two variables and let
$\Pi_n^2$ denote the subspace of polynomials of total degree $n$ in 
$\Pi^2$, which has dimension $\dim \Pi_n^2 = (n+1)(n+2)/2$. If 
$P \in \Pi_n^2$ then 
$$
 P(x,y) = \sum_{k=0}^n \sum_{j=0}^k c_{k,j} x^j y^{k-j}.
$$ 

Let $p$ be a polynomial in one variable and $\theta \in [0,\pi]$. 
For $\xi = (\cos \theta, \sin \theta)$ and $\xb = (x,y)$, the ridge 
polynomial $p(\theta; \cdot)$ is defined by 
$$
   p(\theta;x,y) :=  p (\langle \xb, \xi \rangle) = 
        p(x \cos \theta + y \sin \theta). 
$$
This is a polynomial of two variables and is in $\Pi_n^2$ if $p$ is of 
degree $n$. It is constant on lines that are perpendicular to the 
direction $(\cos \theta,\sin\theta)$.  Especially important to us are the 
polynomials $U_k(\theta; x,y)$, where $U_k$ denotes the Chebyshev 
polynomial of the second kind, 
$$
 U_k (x) = \frac{\sin(k+1)\theta}{\sin \theta}, \qquad   x = \cos \theta.
$$ 
The polynomials $U_k$, $k = 0, 1, 2, \ldots$, are orthogonal polynomials 
with respect to the weight function $\sqrt{1-x^2}$ on $[-1,1]$: 
$$
 \frac{2}{\pi} \int_{-1}^1 U_k(x) U_j(x) \sqrt{1-x^2} dx = \delta_{k,j},
  \qquad k, j \ge 0.
$$
Ridge polynomials play an important role in our study of the Radon 
transform, as can be seen from the following simple fact.

\begin{lem} \label{lem:2.2}  
For $f\in L^1(B^2)$, 
\begin{equation} \label{eq:2.2}
 \int_{B^2} f(x,y) U_k(\phi; x,y)dx dy = \int_{-1}^1 \CR_\phi(f;t) 
      U_k(t) dt. 
\end{equation}
\end{lem} 

\begin{proof}
The change of variables $t = x \cos \phi + y \sin \phi$ and 
$s =-x \sin \phi + y \cos\phi$ amounts to a rotation, which leads to
\begin{align*}
 & \int_{B^2} f(x,y) U_k(\phi; x,y)dx dy  = 
    \int_{B^2} f(t \cos \phi - s \sin \phi, t \sin \phi + s \cos \phi)  
     U_k(t) dt ds \\
 & \qquad  = \int_{-1}^1 \int_{-\sqrt{1-t^2}}^{\sqrt{1-t^2}}
    f(t \cos \phi - s \sin \phi, t \sin \phi + s \cos \phi) ds U_k(t) dt, 
\end{align*}  
the inner integral is exactly $\CR_\phi(f;t)$ by \eqref{eq:Radon}. 
\end{proof}

Let $\CV_k(B^2)$ denote the space of orthogonal polynomials of degree $k$
on $B^2$ with respect to the unit weight function; that is, $P \in \CV_k(B^2)$
if $P$ is of degree $k$ and 
$$
\int_{B^2} P(x,y) Q(x,y) dx dy= 0, \qquad \hbox{for all $Q \in \Pi_{k-1}^2$}.  
$$
A set of polynomials $\{P_{j,k}: 0 \le j \le k\}$ in $\CV_k(B^2)$ is an 
orthonormal basis of $\CV_k(B^2)$ if 
$$
\frac{1}{\pi}\int_{B^2} P_{i,k}(x,y) P_{j,k}(x,y) dx dy = \delta_{i,j}, 
    \qquad  0 \le i,j \le k.
$$
It is known that the polynomials in $\CV_k(B^2)$ are eigenfunctions of a 
second order differential operator $\CD$ (see, for example, \cite{DX}),
$$
  \CD P =  - (k+2) (k-1) P, \qquad \hbox{for all} \quad P\in \CV(B^2);
$$ 
the operator $\CD$ is defined by 
\begin{equation}\label{eq:CD}
\CD: = \Delta^2 - \langle \xb, \nabla\rangle^2 - 2 \langle \xb, \nabla \rangle,
\end{equation}
where $\xb = (x,y)$, $\nabla = (\partial_1, \partial_2)$ is the gradient 
operator with $\partial_1 = \partial / \partial x$  and $\partial_2 = 
\partial / \partial y$, and  $\Delta = \partial_1^2 + \partial_2^2$ is the 
usual Laplacian. It turns out that the Radon projections of orthogonal 
polynomials in $\CV_n(B^2)$ can be computed explicitly. 

\begin{lem} \label{lem:Marr} 
If $P \in \CV_k(B^2)$ then for each $t \in (-1,1)$, $0 \le \theta \le 2\pi$, 
\begin{equation} \label{eq:Marr}
\CR_\theta (P;t) =\frac{2}{k+1} \sqrt{1-t^2} U_k(t) P(\cos \theta,\sin \theta).
\end{equation}
\end{lem} 

\begin{proof}
A change of variable in \eqref{eq:Radon} shows that 
$$
\CR_\theta (P; t ) = \sqrt{1-t^2} \int_{-1}^1
    P \left(t \cos \theta - s \sqrt{1-t^2} \sin \theta, t \sin \theta + 
        s \sqrt{1-t^2} \cos \theta\right) ds.
$$
The integral is a polynomial in $t$ since an odd power of $\sqrt{1-t}$ in the 
integrant is always companied by an odd power of $s$, which has integral zero.
Therefore, $Q(t):= \CR_\theta(P;t) /\sqrt{1-t^2}$ is a polynomial of degree $k$
in $t$ for every $\theta$. Furthermore, the integral also shows that $Q(1) = 
P(\cos \theta, \sin\theta)$. By \eqref{eq:2.2}, 
$$
  \int_{-1}^1  \frac{\CR_\theta(P; t )} {\sqrt{1-t^2}} U_j(t) \sqrt{1-t^2} dt
        =  \int_{B^2} P(x,y) U_j(\theta; x,y) dxdy = 0,
$$
for $ j = 0, 1, \ldots, k-1$, since $P \in \CV_k(B^2)$.  Since $Q$ is of 
degree $k$, we conclude that  $Q (t) = c U_k(t)$ for some constant independent
of $t$. Setting $t =1$ and using the fact that $U_k(1) = k+1$, we have $c =
 P(\cos \theta,\sin\theta)/(k+1)$. 
\end{proof}

This lemma was proved in \cite{Marr}.  It plays an important role in our
development below. We have included its short proof in orde to make the 
paper self-contained.
 
There are several orthogonal or orthonormal bases that are known explicitly
for $\CV_k(B^2)$ (see \cite{DX}). We will work with a special  orthonormal 
basis that is given in terms of the ridge polynomials defined above.  Setting 
$f(x,y) = U_k(\theta;x,y)$ in \eqref{eq:2.2} and using \eqref{eq:Marr}, 
we derive from the orthogonality of the Chebyshev polynomials that
\begin{equation} \label{eq:2.3}
 \frac{1}{\pi} \int_{B^2} U_k(\theta;x,y) U_k(\phi; x,y) dx dy 
           = \frac{1}{k+1} U_k(\cos (\phi - \theta)).  
\end{equation}
Recall that the zeros of $U_k$ are $\cos \theta_{j,k}$, $1 \le j \le k$, where 
$\theta_{j,k} = j \pi /(k+1)$.  As a consequence of \eqref{eq:2.3},  we have 
the 
following result (\cite{LS}, see also \cite{X00}): 

\begin{prop} \label{prop:ortho-basis} 
An orthonormal basis of $\CV_k(B^2)$ is given by
$$
\PP_k: =  \left\{ U_k(\theta_{j,k}; x,y): 0 \le j \le k \right\}, \qquad 
    \theta_{j,k} =\frac{j\pi}{k+1}. 
$$
In particular, the set $\{\PP_k: 0 \le k \le n \}$ is an orthonormal basis 
for $\Pi_n^2$. 
\end{prop}

The polynomials $U_k(\phi;x,y)$ also satisfy a discrete orthogonality. Let 
$$
  \phi_\nu : = \frac{2 \pi \nu}{2m+1}, \qquad\quad  0 \le  \nu \le 2m. 
$$
The discrete orthogonality will follow from the following 
identity:

\begin{prop} \label{prop:sumU}
For $k \ge 0$ and $\theta \in [0,2\pi]$, 
\begin{equation} \label{eq:sumU}
\frac{1}{2m+1} \sum_{\nu = 0}^{2m} U_k(\phi_{\nu}; \cos \theta,\sin\theta)
  U_k(\phi_\nu; x,y) = U_k (\theta; x,y), \qquad x,y \in B^2.
\end{equation}
\end{prop}

\begin{proof}
In order to prove \eqref{eq:sumU} we will need the elementary identities 
\begin{equation}\label{eq:TrigIdent}
 \sum_{\nu=0}^{n} \sin k \phi_\nu = 0 \quad \hbox{and} \quad 
 \sum_{\nu=0}^{n} \cos k \phi_\nu = \begin{cases} 0, &
     \hbox{if $k \ne 0 \mod n+1$} \\
     n+1, &\hbox{if $k = 0 \mod n+1$} \\
   \end{cases}
\end{equation}
that hold for all nonnegative integers $k$. 

Let us denote by $I_k$ the left hand side of \eqref{eq:sumU}. First we 
consider the case of $k = 2l$. Since $U_{2l}$ is an even polynomial, we can
write it as $U_{2l}(t)  = \sum_{j=0}^l b_j t^{2j}$. Using the polar coordinates
$x = r \cos \phi$ and $y = r \sin \phi$, and the fact that 
$(\cos \psi)^{2j}$ can be written as a sum of $\cos 2i \psi$, we can 
rearrange the sum to get 
\begin{align} \label{eq:u1}
  U_{2l}(\theta; x,y) = U_{2l}(r \cos(\theta - \phi)) &  = 
       \sum_{j=0}^l b_j r^{2j} (\cos (\theta - \phi))^{2j} \\ 
  & =  \sum_{j=0}^l b_j(r) \cos 2j (\theta - \phi),  \notag
\end{align}
where $b_j(r)$ is a polynomial of degree $2 j$ in $r$. Furthermore,
we have 
\begin{equation} \label{eq:u2}
 U_{2l}(\phi_\nu;\cos \theta,\sin \theta) = U_{2l}(\cos(\phi_\nu - \theta)) = 
    1 + 2 \sum_{j=1}^l \cos 2j (\theta - \phi_\nu).
\end{equation}
These two expressions allow us to write 
$$
 I_{2l} = \sum_{i= 0}^l d_i \sum_{j =0}^l b_j(r) 
  \frac{1}{2m+1} \sum_{\nu=0}^{2m} \cos 2i (\theta - \phi_\nu)
       \cos 2j (\phi - \phi_\nu)
$$
where $d_0 =1$ and $d_i =2$ for $1 \le i \le l$. The addition formula 
of the cosine function shows that 
\begin{align*}
 & 2 \cos 2i (\theta - \phi_\nu) \cos 2j (\phi - \phi_\nu) \\ 
 & \qquad = \cos 2 ((i+j)\phi_\nu - i\theta - j \phi) + 
      \cos 2 ((i-j)\phi_\nu - i\theta + j \phi) \\
 & \qquad = \cos 2 (i+j)\phi_\nu \cos 2 (i\theta + j \phi) + 
      \sin 2 (i+j)\phi_\nu \sin 2 (i\theta + j \phi) \\ 
 & \qquad \quad + \cos 2 (i-j)\phi_\nu \cos 2(i\theta - j \phi) + 
      \sin 2 (i-j)\phi_\nu \sin 2 (i\theta - j \phi).
\end{align*}
Since $n+1 = 2m+1$ is odd and $\cos 2(i\pm j) \phi_\nu$ have even indices, 
it follows from the above elementary trigonometric identity and 
\eqref{eq:TrigIdent} that 
$$
\frac{1}{2m+1}\sum_{\nu=0}^{2m} \cos 2i (\theta - \phi_\nu) 
    \cos 2j (\phi - \phi_\nu) 
  = \begin{cases} 0, & \hbox{if $i\ne j$}, \\ 
        \frac{1}{2} \cos 2 j(\theta - \phi), & \hbox{if $i = j \ne 0$},\\
         1, & \hbox{if $i = j = 0$.} \end{cases}
$$  
Consequently, using \eqref{eq:u1} again, we conclude that 
$$
  I_{2l} = \sum_{j=0}^l b_j(r) \cos 2 j(\theta - \phi) =
     U_{2l}(r \cos (\theta - \phi)) = U_{2l}(\theta; x,y).
$$
This completes the proof for the case $k = 2l$. For the case $k = 2l-1$,
we need to use the identities 
$$
 U_{2l-1} (\phi_\nu;x,y) =  U_{2l}(r \cos(\phi_\nu - \phi))  
    =  \sum_{j=1}^l b_j(r) \cos (2j-1)(\phi_\nu - \phi),  
$$
where $b_j(r)$ is a polynomial of degree $2 j-1$ in $r$, which is 
derived using the fact that $(\cos \theta)^{2j-1}$ can be written as a sum
of $\cos (2 i -1) \theta$. Furthermore, we have 
$$
U_{2l-1}(\phi_\nu;\cos \theta,\sin \theta) = 
  U_{2l-1}(\cos(\phi_\nu - \theta)) = 
    2 \sum_{j=1}^l \cos (2j-1) (\theta - \phi_\nu).
$$
The rest of the proof will follow as in the case $k =2l$.
\end{proof}

\begin{cor} \label{cor:discretU}
For $k \ge 0$ and $0 \le i,j\le k$,
\begin{equation} \label{eq:discretU}
\frac{1}{2m+1} \sum_{\nu = 0}^{2m} 
  U_k(\phi_{\nu}; \cos \theta_{i,k},\sin\theta_{i,k})
  U_k(\phi_{\nu}; \cos \theta_{j,k},\sin\theta_{j,k}) = (k+1)\delta_{i,j}.
\end{equation}
\end{cor}

\begin{proof}
This follows from setting $\theta = \theta_{i,k}$, $x = \cos\theta_{j,k}$
and $y = \sin \theta_{j,k}$ in the identity \eqref{eq:sumU}, and using the
fact that 
$$
  U_k(\theta_{i,k}; \cos\theta_{j,k}, \sin\theta_{j,k}) = 
   U_k (\cos \theta_{i-j,k}) = (k+1) \delta_{i,j},
$$ 
where in the last step we have used the fact that $U_k(0) = (k+1)$.
\end{proof}

\begin{rem} \label{rem:2.1}
The identity \eqref{eq:sumU} is established for the case $n =2m$. It should
be pointed out that the result does not hold for $n = 2m - 1$. In fact, if
$n = 2m-1$, then following the above proof leads to an identity analogous 
to \eqref{eq:sumU} only for $0 \le k \le m-1$. More precisely, we have 
$$
  \frac{1}{2m} \sum_{\nu = 0}^{2m -1}  U_k(\phi_{\nu}; \cos \theta,\sin\theta)
  U_k(\phi_\nu; x,y) = U_k (\theta; x,y), \qquad x,y \in B^2,
$$ 
holds for $0 \le k \le m-1$. Indeed, in the case $m = 2 p$ it can be proved
that 
$$
\frac{1}{2m}\sum_{\nu=0}^{2m-1} U_m (\phi_{\nu}; \cos \theta,\sin\theta)
 U_m (\phi_{\nu}; \cos \phi,\sin\phi) = 
    U_m (\cos(\theta-\phi)) + 2 T_m (\cos (\theta - \phi)), 
$$
where $T_m(x) = \cos m \psi$, $x = \cos \psi$, is the Chebyshev polynomial 
of the first kind. 
\end{rem}

There is yet another remarkable identity for $U_k(\cdot;x,y)$, 
proved in \cite{X00}, which takes the summation over $\theta_{j,k}$ instead
of taking over $\phi_\nu$: 

\begin{lem} \label{lem:2.6}
For $k \ge 0$, 
\begin{equation}\label{eq:sumU2}
\sum_{j=0}^k U_k(\theta_{j,k}; x,y) 
             U_k(\theta_{j,k}; \cos \phi,\sin\phi) 
        = (k+1) U_k(\phi;x,y).
\end{equation}
\end{lem}

This identity is a consequence of the compact formula for the 
reproducing kernel of $\CV_n(B^2)$ in \cite{X99}. It will also play an 
important role in our development below.

The two remarkable identities \eqref{eq:sumU} and \eqref{eq:sumU2} 
hold the key for our new algorithms. The orthogonal polynomials in 
$\CV_n(B^2)$ are usually studied as a special case of the orthogonal 
polynomials with respect to the weight functions $W_\mu(x,y) : = 
(1-x^2-y^2)^{\mu-1/2}$, $\mu \ge 0$. Most of their properties are shared
by orthogonal families associated with $W_\mu$ for all $\mu \ge 0$,
and furthermore, many properties can be extended to higher dimensions
(see \cite{DX}). The orthogonality in \eqref{eq:discretU}, however, is very 
special; it is not shared by any other orthogonal families associated 
with $W_\mu$ for $\mu \ne 1/2$ and there is no direct extension to 
higher dimensions; see, for example, the discussion in  \cite{X00}.

\section{Fourier Orthogonal expansions}
\setcounter{equation}{0}

\subsection{ Orthogonal expansions on the disk}  \label{sec:S-B}

The standard Hilbert space theory shows that any function in $L^2(B^2)$ can be 
expanded as a Fourier orthogonal series in terms of $\CV_n$. More precisely, 
\begin{equation} \label{eq:fourier}
   L^2(B^2) = \sum_{k=1}^\infty \bigoplus \CV_k(B^2): 
\qquad   f = \sum_{k=1}^\infty \proj_k f, 
\end{equation} 
where $\proj_k f$ is the orthogonal projection of $f$ onto the subspace 
$\CV_k(B^2)$. Our reconstruction algorithm will be based on this Fourier 
orthogonal expansion. It is well known that $\proj_n f$ can be written as 
an integral operator in terms of 
the reproducing kernel of $\CV_k(B^2)$ in $L^2(B^2)$. A compact formula
of the reproducing kernel is given in \cite{X99}. For our purpose, we seek 
a formula for $\proj_n f$ that will relate it to the Radon transform. 

In terms of the special orthonormal basis $\{\PP_k: k \ge 0\}$ given in 
Proposition \ref{prop:ortho-basis}, we can write $\proj_k f$ as  
\begin{equation} \label{eq:fourierProj}
  \proj_k f =  \sum_{j=0}^k  \wh f_{j,k} U_k(\theta_{j,k}; \cdot), 
 \qquad \wh f_{j,k} = \frac{1}{\pi} \int_{B^2} f(x,y)U_k(\theta_{j,k}; x,y) 
 dx\, dy.
\end{equation} 
The remarkable identity \eqref{eq:sumU} allows us to express
the Fourier coefficients in terms of Radon projections:

\begin{prop} \label{prop:fcoeff}
Let $m$ be a nonnegative integer. For $0 \le k \le 2m$, the Fourier 
coefficient $\wh f_{j,k}$ satisfies 
\begin{equation} \label{eq:fcoeff}
\wh f_{j,k} = \frac{1}{2m+1} \sum_{\nu =0}^{2m}
 \frac{1}{\pi} \int_{-1}^1 \CR_{\phi_\nu} (f;t) U_k(t) dt \, 
   U_k(\cos (\theta_{j,k} - \phi_\nu)).
\end{equation}
\end{prop}

\begin{proof}
Let us denote the right hand side of \eqref{eq:fcoeff} by $g_{j,k}$. By
the equation \eqref{eq:2.2}, 
\begin{align*}
g_{j,k} & =  \frac{1}{2m+1} \sum_{\nu =0}^{2m}
             \frac{1}{\pi} \int_{B^2} f(x,y) U_k(\phi_\nu;x,y) dx\, dy \, 
             U_k(\cos (\theta_{j,k} - \phi_\nu)) \\  
     & =  \frac{1}{\pi} \int_{B^2} f(x,y)
         \left[ \frac{1}{2m+1} \sum_{\nu =0}^{2m} U_k(\phi_\nu;x,y) 
         U_k(\cos (\theta_{j,k} - \phi_\nu)) \right]dx \,dy   \\
     & =  \frac{1}{\pi} \int_{B^2} f(x,y) U_k(\theta_{j,k}; x,y)dx\,dy
\end{align*}
by the identity \eqref{eq:sumU}. Hence, by the definition of $\wh f_{j,k}$, 
$g_{j,k} = \wh f_{j,k}$. 
\end{proof}

A further application of the identity \eqref{eq:sumU2} leads us to  
the expression of the Fourier projection operator in terms of 
Radon projections:

\begin{thm}  \label{thm:proj-sum}
For $m \ge 0$ and $k \le 2m$,  the operator $\proj_k f$ can be 
written as
\begin{equation} \label{eq:proj} 
\proj_{k} f(x,y) = \frac{1}{2m+1} \sum_{\nu=0}^{2m} \frac{1}{\pi} 
  \int_{-1}^1 \CR_{\phi_\nu}(f;t) U_k(t) dt (k+1) U_k(\phi_\nu;x,y) .
\end{equation}  
\end{thm}

\begin{proof}
By \eqref{eq:fcoeff} and the identity \eqref{eq:sumU2} in Lemma \ref{lem:2.6}, 
\begin{align*}
\proj_k f(x,y) = & 
   \frac{1}{2m+1}  \sum_{\nu = 0}^{2m} \frac{1}{\pi} \int_{-1}^1 
      \CR_{\phi_\nu} (f;t) U_k(t) dt \\
    &\qquad\qquad  \times \sum_{j=0}^k  
         U_k(\cos (\theta_{j,k} - \phi_\nu)) U_k(\theta_{j,k};x,y) \\
  =  & \frac{1}{2m+1}  \sum_{\nu = 0}^{2m} \frac{1}{\pi} \int_{-1}^1 
      \CR_{\phi_\nu} (f;t) U_k(t) dt (k+1) U(\phi_\nu;x,y),  
\end{align*}
 which is what we need to prove. 
\end{proof}

We denote the $n$-th partial sum of the expansion \eqref{eq:fourier} 
by $S_n f$; that is,
$$
S_n f(x,y) =  \sum_{k=0}^n \proj_k f (x,y).
$$
The operator $S_n$ is a projection operator from $L^2(B^2)$ onto $\Pi_n^2$. 
An immediate consequence of the Theorem \ref{thm:proj-sum} is the following
result that will play an essential role in dering the new algorithm:

\begin{cor}  \label{thm:partial-sum}
For $m \ge 0$, the partial sum operator $S_{2m} f$ can be written as
\begin{equation} \label{eq:S2m} 
S_{2m}(f;x,y) =  \sum_{\nu=0}^{2m} \frac{1}{\pi} 
  \int_{-1}^1 \CR_{\phi_\nu}(f;t) \Phi_\nu (t; x,y) dt
\end{equation}  
where 
\begin{equation} \label{eq:Phi} 
 \Phi_{\nu}(t; x,y) = \frac{1}{2m+1} \sum_{k=0}^{2m} (k+1) 
     U_k(t)U_k(\phi_\nu;x,y).
\end{equation} 
\end{cor}
 
The identity \eqref{eq:S2m} shows that $S_{2m}f$ can be expressed in 
terms of Radon projections in $2m+1$ directions. This result, previously
unnoticed, holds the key for our new algorithm given in Section 4. It should 
be pointed out (recall Remark \ref{rem:2.1}) that such an identity does not 
hold for $S_{2m-1} f$. 

\subsection{ Summability of orthogonal expansions}  \label{sec:S-eta}
Let $L^p(B^2)$ denote the usual $L^p$ space on $B^2$ with norm 
$\|\cdot\|_p$ for $1 \le p < \infty$ and identify it with $C(B^2)$ of
continuous functions with uniform norm for $p = \infty$. If $f\in 
L^p(B^2)$, the error of best approximation by 
polynomials of degree at most $n$ is defined by 
\begin{equation}\label{eq:bestapp}
   E_n(f)_p : = \inf \left \{ \|f - P \|_p:  P \in \Pi_n^2 \right\}.
\end{equation}
It is well-known that $E_n(f)_p \to 0$ for $f \in L^p(B^2)$ as $n \to \infty$. 
The partial sum $S_n f$ of the orthogonal expansion is the best 
approximation to $f$ in the $L^2$ norm; that is, 
$$
      \|f - S_n f\|_2 = E_n(f)_2, \qquad f \in L^2(B^2).
$$  
However, the partial sum $S_n f$ does not converge to $f$ pointwisely 
if $f$ is merely continuous; see \cite{X01} and the discussion below 
in Section 5. 

To study the convergence of our orthogonal expansions we will introduce 
some summability methods. Such a method takes the form 
$$
  \sum_{j=0}^\infty  a_{j,n} S_j (f), \qquad a_{j,n} \in \RR
   \quad\hbox{and}\quad    \sum_{j=0}^\infty a_{j,n}=1.
$$
Many summability methods, for example, the Poisson means and the 
Ces\`aro means,  have better convergence behavior (see \cite{X99}).  
For our purpose we will use methods for which the operators are polynomials 
and we would still want to retain the property that polynomials up to certain 
degree are preserved. 

It turns out that this can be done quite easilly using a multiplier function. 

\begin{defn}
Let $r$ be a positive integer, and let $\eta \in C^r[0, \infty)$. Then
$\eta$ is called a {\it multiplier function} if 
$$
  \eta(t) = 1, \quad 0 \le t \le 1, \quad\hbox{and}\quad
  \supp \eta \subset [0,2].  
$$
\end{defn} 

Let $\eta$ be a multiplier function. We define an operator $S_{2m}^\eta$ by
$$
 S_{2m}^\eta (f; x,y) = \sum_{k=0}^{2m} \eta\left(\frac{k}{m}\right)
      \proj_k f(x,y).     
$$
Such an operator was used in the literature for approximation by 
spherical polynomials on the unit sphere, and it was been used for various 
other orthogonal expansions in \cite{X05}, including the expansions
on the unit disk.  The following theorem is essentially contained in 
\cite{X05}. 

\begin{prop} \label{prop:S-eta}
Let $\eta \in C^3[0,\infty)$ be a multiplier function. Let $f \in L^p(B^2)$, 
$1 \le p \le \infty$. Then  
\begin{enumerate}
\item $S_{2m}^\eta f \in \Pi_{2m}^2$ and $S_{2m}^\eta P  = P$ for $P\in 
\Pi_m^2$;
\item For $m \ge 0$ there is a constant $c$ such that 
$$
 \|S_{2m}^\eta f \|_p \le c \|f\|_p \quad \hbox{and} \quad  
 \|f- S_{2m}^\eta f \|_p \le c E_m (f)_p.
$$
\end{enumerate}
\end{prop}

This means that the operator $S_{2m}^\eta$ is very well behaved: it 
preserves polynomials of degree up to $m$ and it approximates $f$ as 
accurate as any polynomial of degree at most $m$. Using
Theorem \ref{thm:proj-sum} we also have the following:

\begin{thm}  \label{thm:S-eta}
For $m \ge 0$, the operator $S_{2m}^\eta$ can be written as
\begin{equation} \label{eq:S2m-eta} 
S_{2m}^\eta(f;x,y) =  \sum_{\nu=0}^{2m} \frac{1}{\pi} 
  \int_{-1}^1 \CR_{\phi_\nu}(f;t) \Phi_\nu^\eta (t; x,y) dt
\end{equation}  
where 
\begin{equation} \label{eq:Phi-eta} 
 \Phi_{\nu}^\eta (t; x,y) = \frac{1}{2m+1} \sum_{k=0}^{2m} 
      \eta\left(\frac{k}{m}\right) (k+1)U_k(t_j)U_k(\phi_\nu;x,y).
\end{equation} 
\end{thm}

\begin{rem} \label{rem:3.1}
We can also use other summability methods, not necessarily prescribed 
by the multiplier function. The essence of Proposition \ref{prop:S-eta},
however, should be preserved. 
\end{rem}

\subsection{Fourier orthogonal expansion on a cylinder domain}\label{sec:S-BL}
                    
Let $L > 0$ and let $B_L$ be the cylinder region 
$$
B_L := B^2 \times [0,L] = \{(x,y,z): (x,y) \in B^2, \,\, z\in [0,L]\}. 
$$ 
Using the result in the Subsection \ref{sec:S-B} we can also get an expression 
for the partial sum operator on $B_L$, which will lead us to a 3D 
reconstruction
algorithm. We consider orthogonal polynomials with respect to the inner
product
\begin{equation} \label{eq:inner}
 \langle f, g\rangle_{B_L}  =  \frac{1}{\pi} \int_{B_L} f(x,y,z) g(x,y,z)
       W_L(z)\,dx\, dy \,dz,
\end{equation}
where $W_L$ is a nonnegative function defined on $[0,L]$ with all its 
moments on $[0,L]$ assumed finite and normalized so that  
$\int_0^L W_L(z) dz =1$.

Let $\Pi_n^3$ denote the space of polynomials of total degree at most $n$ 
in three variables.  Let $\CV_n(B_L)$ denote the subspace of orthogonal 
polynomials of degree $n$ on $B_L$ with respect to the inner product 
\eqref{eq:inner};  
that is, $P \in \CV_n(B_L)$ if $\langle P, Q\rangle_{B_L} =0$ for all 
polynomial $Q \in \Pi_{n-1}^3$. 

Let $p_k$ be the orthonormal polynomials with respect to $W_L$ on 
$[0,L]$. Let $U_k(\theta_{j,k};x,y)$ be defined as in the previous subsection. 

\begin{prop}
An orthonormal basis for $\CV_n(B_L)$ is given by 
$$
\PP_n:=\{P_{n,k,j}: 0 \le j \le k \le n\}, \qquad 
    P_{n,k,j}(x,y,z) = p_{n-k}(z) U_k(\theta_{j,k};x,y).   
$$
In particular, the set $\{\PP_l: 0 \le l \le n \}$ is an orthonormal basis 
for $\Pi_n^3$. 
\end{prop}

This is an easy consequence of the fact that $B_L$ is a product of $B^2$
and $[0,L]$. For $f \in L^2(B_L)$, the Fourier coefficients of $f$ with 
respect to the orthonormal system $\{\PP_l: l \ge 0\}$ are given by 
$$
  \wh f_{l,k,j} =\frac{1}{\pi} \int_{B_L} f(x,y,z)P_{l,k,j}(x,y,z) dxdy 
       W_L(z) dz,  \quad 0 \le j \le k \le  l.
$$ 
Let $S_n f$ denote the Fourier partial sum operator, 
$$
S_n f(x,y,z) = \sum_{l=0}^n \sum_{k=0}^l \sum_{j=0}^k \wh f_{l,k,j}
      P_{l,k,j}(x,y,z). 
$$
Just like its counterpart in two variables,  this is a projection operator 
and is independent of the particular choice of the bases of $\CV_n(B_L)$. 

We retain the notation  $\CR_\phi(g;t)$ for the Radon projection of a 
function $g: B^2 \mapsto \RR$. For a fixed $z$ in $[0,L]$, we define
\begin{equation}\label{eq:Radon3D}
  \CR_\phi(f(\cdot,\cdot,z); t) := \int_{I(\phi,t)} f(x,y,z) dx dy,
\end{equation}
which is the Radon projection of $f$ in a disk that is perpendicular to 
the $z$-axis.  

The following is an analogue of Theorem \ref{thm:partial-sum} for the
cylinder $B_L$.

\begin{thm} \label{thm:partial-sum3D}
For $m \ge 0$, 
\begin{align} \label{eq:S2m3D}
 & S_{2m} f(x,y,z) \\ &
  \qquad = \frac{1}{2m+1} \sum_{\nu=0}^{2m} 
    \frac{1}{\pi} \int_{-1}^1 \int_0^L \CR_{\phi_\nu} (f(\cdot,\cdot,w);t)
      \Phi_\nu (w,t;x,y,z)  W_L(w)dw\,dt \notag
\end{align}
where
\begin{equation} \label{eq:Phi3D}
  \Phi_\nu(w,t;x,y,z) = \sum_{k=0}^{2m}  (k+1) U_k(t) U_k(\phi_\nu;x,y) 
      \sum_{l=0}^{2m-k} p_l(w) p_l(z).
\end{equation}
\end{thm} 

\begin{proof}
Using Proposition \ref{prop:fcoeff} and the product nature of the region, 
$$
 \wh f_{l,k,j} =  \frac{1}{2m+1} \sum_{\nu =0}^{2m}
   \frac{1}{\pi} \int_{-1}^1 \CR_{\phi_\nu} (f_{l-k};t) U_k(t) dt \, 
     U_k(\cos (\theta_{j,k} - \phi_\nu) ),
$$
where
$$
   f_l(x,y) = \frac{1}{\pi} \int_0^L f(x,y,w) p_l(w)  W_L(w) dw, 
     \qquad  l \ge 0. 
$$
Substituting this expression of $\wh f_{l,k,j}$ into the formula of
$S_{2m} f$ and using the identity \eqref{eq:sumU2} in Lemma \ref{lem:2.6}, 
we then obtain 
\begin{align*}
&  S_{2m} f(x,y,z) =  \frac{1}{2m+1} \sum_{\nu =0}^{2m}
   \frac{1}{\pi} \int_{-1}^1 \sum_{l=0}^{2m} \sum_{k=0}^l
     \CR_{\phi_\nu} (f_{l-k};t) U_k(t) dt \\ 
      & \qquad\quad\qquad \qquad\quad\qquad \qquad\quad
\times (k+1) U_k(\phi_\nu;x,y) p_{l-k} (z)\\
&\qquad = \frac{1}{2m+1} \sum_{\nu =0}^{2m}
     \frac{1}{\pi} \int_{-1}^1 \sum_{k=0}^{2m} \sum_{l=0}^{2m-k} 
        \CR_{\phi_\nu} (f_{l};t)  p _{l} (z) 
   (k+1) U_k(\phi_\nu;x,y) U_k(t) dt,
\end{align*} 
where in the second equality we have exchanged the two inner summations. 
From \eqref{eq:Radon3D} and the definition of $f_l$ it is easy to see that 
$$
 \CR_{\phi_\nu} (f_l;t) = \frac{1}{\pi} \int_0^L 
    \CR_{\phi_\nu} (f(\cdot,\cdot, w); t) p_l(w) W_L(w)\,dw,
$$
from which the proof follows upon rearranging terms in the summation. 
\end{proof}

\begin{rem}
Clearly one can also consider summability of orthogonal expansions on $B_L$. 
For example, one can define the operator with multiplier factors just as in the
case of orthogonal expansion on $B^2$. We shall not elaborate.
\end{rem}

\section{New Reconstruction algorithms}
\setcounter{equation}{0}

\subsection{ Reconstruction Algorithm for 2D images} 
The identity \eqref{eq:S2m} expresses $S_{2m} f$ in terms of the Radon 
projections $\CR_{\phi_\nu} (f;t)$ of $2m+1$ equally spaced angles 
$\phi_\nu$, $0\le \nu \le 2m$, along the circumference of the disk. These 
Radon projections are defined for all parameters $t$. In order to make use 
of the Radon data from the parallel geometry,  we will use a quadrature 
rule to get a discrete approximation of the integrals 
$$
\int_{-1}^1  \CR_{\phi_\nu} (f;t) U_k(t) dt = \int_{-1}^1 
    \frac{\CR_{\phi_\nu} (f;t)}{\sqrt{1-t^2} } U_k(t) \sqrt{1-t^2}  dt 
$$
in \eqref{eq:S2m}. The result will be our algorithm.

If $f$ is a polynomial then the equation  \eqref{eq:Marr} shows that 
$\CR_\phi (f; t)/\sqrt{1-t^2}$  is also a polynomial. Hence,  we choose 
a quadrature rule for the integral with respect to $\sqrt{1-x^2}$ on $[-1,1]$. 
Let us denote such a quadrature rule by $\CI_n g$; then 
\begin{equation} \label{eq:quadrature}
\frac{2}{\pi} \int_{-1}^1 g(t) \sqrt{1-t^2}dt \approx 
   \sum_{j=1}^n \lambda_j g(t_j) := \CI_n(g),   
\end{equation}
where $t_1, \ldots, t_n$ are distinct points in $(-1,1)$ and $\lambda_j$
are real numbers such that $\sum_{j=1}^n \lambda_j =1$. If equality holds
in \eqref{eq:quadrature} whenever $g$ is a polynomial of degree at most 
$\rho$, then the quadrature rule is said to have the degree of exactness
$\rho$. 

Among all quadrature rules with a fixed number of nodes, the Gaussian 
quadrature has the highest degree of exactness. It is given by
\begin{equation} \label{eq:gaussian}
 \frac{1}{\pi} \int_{-1}^1 g(t) \sqrt{1-t^2} dt = \frac{1}{n+1}
    \sum_{j=1}^{n} \sin^2 \frac{j\pi}{n+1} 
       g\left(\cos \frac{j\pi}{n+1}\right) : = \CI_n^G(g)
\end{equation}                                    
for all polynomials $g$ of degree at most $2n-1$; that is, its degree of 
exactness is $2n-1$.  Note that $j \pi /(n+1)$ are zeros of the Chebyshev
polynomial $U_n$. 

Using quadrature formula in \eqref{eq:S2m} gives our  reconstruction
algorithm, which produces a polynomial $\CA_{2m}f$ defined below.

\begin{algo}  \label{algo:2D} 
Let the quadrature rule be given by \eqref{eq:quadrature}. For $m \ge 0$ and
$ (x,y) \in B^2$,
\begin{equation} \label{eq:Algo} 
\CA_{2m}(f;x,y) =  \sum_{\nu=0}^{2m} \sum_{j=1}^{n}
       \CR_{\phi_\nu}(f;t_j) T_{j,\nu}(x,y), 
\end{equation}  
where 
$$
 T_{j,\nu}(x,y) = \frac{\lambda_j}{2 (2m+1) \sqrt{1-t_j^2}}\Phi_\nu (t_j;x,y).
$$
\end{algo}

For a given $f$, the approximation process $\CA_{2m} f$ uses the Radon data 
$$
  \{\CR_{\phi_\nu} (f; t_j): 0 \le \nu \le 2m, \quad
         1 \le j \le n\}
$$
of $f$. The data consists of Radon projections on $2m+1$ equally spaced 
directions along the circumference of the disk (specified by $\phi_\nu$)
and there are $n$ parallel lines (specified by $t_j$) in each 
direction. If these parallel Radon projections are taken from an image
$f$, then the algorithm produces a polynomial $\CA_{2m} f$ which gives an 
approximation to the original image. 

The polynomial $\CA_{2m}$ is particularly handy for numerical implementation, 
since one could save $T_{j,\nu}$ in a table beforehand. This provides 
a very simple algorithm: given the Radon data, one only has to perform
addition and multiplication to evaluate $\CA_{2m}$ in \eqref{eq:Algo} to get 
a reconstruction of image.
  
A good choice of the quadrature rule is Gaussian quadrature. We choose in 
particular $n=2m$ so that the nodes of the quadrature rule \eqref{eq:gaussian} 
becomes $t_j = \cos \theta_{j,2m} = \cos j \pi/(2m+1)$. In this case, our 
algorithm takes a particular simple form. 

\begin{algo} \label{algo:2D-G} 
For $m \ge 0$, $(x,y) \in B^2$,
\begin{equation} \label{eq:AlgoG} 
 \CA_{2m}(f; x,y) = \sum_{\nu=0}^{2m} \sum_{j=1}^{2m}
   \CR_{\phi_\nu}(f;\cos \theta_{j,2m}) T_{j,\nu}(x,y),
\end{equation} 
where 
\begin{equation} \label{eq:TjvG} 
 T_{j,\nu}(x,y) = \frac{1}{(2m+1)^2} 
     \sum_{k=0}^{2m} (k+1)\sin((k+1)\theta_{j,2m})U_k(\phi_\nu;x,y).
\end{equation}
\end{algo}

\begin{rem}
It is tempting to choose $n=2m+1$ in the Gaussian quadrature. In fact, 
such a choice has one advantage: the operator $\CA_{2m}$ will be 
a projection operator onto $\Pi_{2m}^2$.  We choose $n =2m$ so that
$\phi_\nu = \frac{2\nu \pi}{2m+1}$ and $\theta_{j,2m} = \frac{j \pi}{2m+1}$ 
have common denominator.  It also turns out that this choice works 
perfectly with the fan beam geometry of the projections, which will be 
reported elsewhere. 
\end{rem}

The convergence of the algorithm will be discussed in the following 
section. Here we state one result that is a simple consequence of the
definition.

\begin{thm} \label{thm:A2m}
The operator $\CA_{2m}$ in Algorithm \ref{algo:2D} preserves polynomials
of degree $\sigma$. More precisely, $\CA_{2m} (f) = f$ whenever $f$ is a
polynomial of degree at most $\sigma$. In particular, the operator 
$\CA_{2m}$ in Algorithm \ref{algo:2D-G} preserves polynomials of degree 
at most $2m-1$. 
\end{thm}

\begin{proof}
The polynomial $\CA_{2m} f$ is obtained by applying quadrature rule
\eqref{eq:quadrature} to \eqref{eq:S2m}. If $f$ is a polynomial of degree
at most $\sigma$, then so is $\CR_{\phi_\nu} (f;t)/\sqrt{1-t^2}$ for every 
$\nu$. Furthermore, the polynomial $\Phi_\nu(t;x,y)$ is a polynomial of 
degree $2m$ in $t$. Hence, when we apply the quadrature rule
\eqref{eq:quadrature}
to the integral \eqref{eq:S2m}, the result is exact. Therefore, $\CA_{2m} f = 
S_{2m}f = f$. 
\end{proof}

If the quadrature rule is the Gaussian quadrature in \eqref{eq:gaussian}, then 
the polynomials preserved by the operator $\CA_{2m}$ have the highest 
degrees among all quadrature rules that use the same number of nodes. Such 
a choice will ensure better approximation behavior of $\CA_{2m}$.

\begin{rem} 
Using the angles $\phi_\nu$, another projection operator, call it $\CJ_{2m}f$, 
has been constructed in \cite{BX} based on the parallel geometry. For almost
all choices of $\{t_j \in (-1,1): 0 \le j \le m\}$, the operator $\CJ_{2m}$ is 
the unique polynomial of degree $2 m$ determined by the conditions 
$$
  \CR_{\phi_\mu} (\CJ_{2m} f; t_i) =   \CR_{\phi_\mu} (f; t_i), \qquad 
         0 \le \nu \le 2m, \quad 0 \le j \le m.
$$
However, the construction of $\CI_{2m}$ requires solving a family of linear
system of equations whose coefficient matrices, depending on the choice of 
$t_j$,  appear to be badly ill-conditioned.Ê
\end{rem}

\subsection{ Reconstruction algorithm for 2D images with a multiplier function}

Instead of using \eqref{eq:S2m}, we can also start from the summabilty
with the multiplier factor \eqref{eq:S2m-eta}.  The result is another 
reconstruction algorithm. Here we state the resulting algorithm only 
for the Gaussian quadrature \eqref{eq:gaussian} with $n = 2m$.

\begin{algo} \label{algo:2D-eta} 
For $m \ge 0$, $(x,y) \in B^2$,
\begin{equation} \label{eq:Algo -eta} 
 \CA_{2m}^\eta(f; x,y) = \sum_{\nu=0}^{2m} \sum_{j=1}^{2m}
     \CR_{\phi_\nu}(f;\cos \theta_{j,2m}) T_{j,\nu}^\eta(x,y),
\end{equation} 
where 
$$
 T_{j,\nu}^\eta(x,y) = \frac{1}{(2m+1)^2} 
     \sum_{k=0}^{2m}\eta\left(\frac{k}{m}\right)(k+1)
     \sin((k+1)\theta_{j,2m})U_k(\phi_\nu;x,y).
$$
\end{algo}

For a given $f$, the approximation process $\CA_{2m}^\eta f$ uses the same 
Radon data of $f$ as $\CA_{2m}f $. It also has the same simple structure for
numerical implementation. Its approximation behavior appears to be better
than that of $\CA_{2m} f$. We conclude this subsection with the following 
analogous of Theorem \ref{thm:A2m}:

\begin{thm} \label{thm:A2m-eta}
The operator $\CA_{2m}^\eta$ preserves polynomials of degree $m$. More 
precisely, $\CA_{2m}^\eta (f) = f$ whenever $f$ is a polynomial of degree at 
most $m$. 
\end{thm}

We can obtain other reconstruction algorithms using different summability
methods; see, however, Remark \ref{rem:3.1}.

\subsection{ Reconstruction algorithm for 3D images} 

In order to get a reconstruction algorithm for 3D images on the cylinder 
region $B_L$, we choose the weight function to be 
$$
      W_L (z) = \frac{1}{\pi} \frac{1}{\sqrt{z(L-z)}}, \qquad z \in [0,L],
$$ 
which is the Chebyshev  weight function on the interval $[0,L]$, normalized
so that  its integral over $[0,L]$ is 1. 
Let $T_k$ be the Chebyshev polynomial of the first kind. Define $\wt T_k$ by
$$
  \wt T_0(z) =1, \qquad \wt T_k(z) = \sqrt{2} T_k(2 z/L -1), \quad k \ge 1.
$$
The polynomials $\wt T_k$ are orthonormal polynomials with respect to 
$W_L$ on $[0,L]$. 

 To obtain an algorithm using parallel Radon projections, we start from 
 \eqref{eq:S2m3D} and apply quadrature rules on its integrals.  
For the integral in $z$, we use the Gaussian quadrature for $W_L$. Set 
$$
  \xi_{i,n} = \frac{(2i+1) \pi}{2n} \quad \hbox{and} \quad
   z_i = \frac{1+\cos\xi_{i,n}}{2}, \qquad 0 \le i \le n-1, 
$$
where $z_i$ are zeros of $T_n(z)$; then the Gaussian quadrature on $[0,L]$
takes the form,
\begin{equation} \label{eq:gaussianT}
  \frac{1}{\pi} \int_{0}^L g(z)\frac{dz} {\sqrt{z(L-z)}} = 
      \frac{1}{n} \sum_{i=0}^{n-1} g (z_i),
\end{equation} 
which holds whenever $g$ is a polynomial of degree at most $2n-1$. For 
the integral in $t$, we could use the quadrature rule \eqref{eq:quadrature}. 
For simplicity, however, we will only use the Gaussian quadrture 
\eqref{eq:quadrature}.  

This way we get an algorithm for reconstruction of images on $B_L$. 
The algorithm produces a polynomial $\CB_{2m}$ of three variables as follows: 

\begin{algo} \label{algo:3D} 
Let $\gamma_{\nu,j,i}=\CR_{\phi_\nu} (f(\cdot,\cdot, z_i);\cos \theta_{j,2m})$.
For $m \ge 0$, 
$$
 \CB_{2m} f(x,y,z) := \sum_{\nu=0}^{2m} \sum_{j=1}^{2m} 
   \sum_{i=0}^{n-1} \gamma_{\nu,j,i}  T_{\nu,j,i} (x,y,z),        
$$
where 
$$
T_{\nu,j,i}(x,y,z) = \frac{1}{n(2m+1)} \Phi_\nu(z_i, \cos\theta_{j,2m};x,y,z). 
$$
\end{algo}

For a given function $f$, the approximation process $\CB_{2m}$ uses the Radon 
data 
$$
\{\CR_{\phi_\nu} (f(\cdot,\cdot,z_i;\cos\theta_{j,2m}): 0 \le\nu\le 2m, \quad
         1 \le j \le 2m, \quad 0 \le i \le n-1 \}
$$
of $f$. The data consists of Radon projections on $n$ disks that are
perpendicular to the $z$-axis (specified by $z_i$); on each disk the Radon
projections are taken in $2m+1$ equally spaced directions along the 
circumference of the disk (specified by $\phi_\nu$) and $2m$ parallel 
lines (specified by $\cos \psi_j$) in each direction. 
We can use this approximation for the reconstruction of the 3D images from
parallel X-ray data. In practice, the integer $n$ of $z$ direction should
be chosen so that the resolution in the $z$ direction is comparable to the 
resolution on each disk to achieve isotropic result. 

The following theorem is an analogous of Theorem \ref{thm:A2m} for $B_L$.

\begin{thm} \label{thm:B2m}
If $n \ge 2m$, then the operator $\CB_{2m}$ in Algorithm \ref{algo:3D} 
preserves polynomials of degree $2m-1$. More precisely, $\CB_{2m} (f) = f$ 
whenever $f$ is a polynomial of degree at most $2m-1$. 
\end{thm} 

\begin{proof}
If $f$ is a polynomial of degree at most $2m-1$, then 
$\CR_{\phi}(f(\cdot,\cdot, w);t)/\sqrt{1-t^2}$ is a polynomial of degree
at most $2m-1$ both in the $t$ variable and in the $w$ variable. The function
$\Phi_\nu(w,t;x,y,z)$ is of degree $2m$ in both the $t$ variable and the $w$
variable. Hence, when we use the quadrature rules for $S_{2m}f$ in 
\eqref{eq:S2m3D}, the result is exact if the quadrature rules are exact 
for polynomials of degree $4m-1$. For the quadrature rule
\eqref{eq:gaussianT} this holds if $n \ge 2m$.  
\end{proof} 

\begin{rem} 
In the $z$ direction, we choose the weight function 
$(z(1-z))^{-1/2}$ instead of constant weight function. The reason
lies in the fact that the Chebyshev polynomials of the first kind are
simple to work with and the corresponding Gaussian quadrature 
\eqref{eq:gaussianT} is explicit. If we were to use constant weight 
functions, we would have to work with Legendral polynomials, whose zeros
(the nodes of Gaussian quadrature) can be given only numerically.  
\end{rem}

\section{Convergence of the reconstruction algorithm}
\setcounter{equation}{0}

In this section we study the convergence behavior of the reconstruction 
algorithm for 2D images and we work with $\CA_{2m}$ in \eqref{eq:AlgoG},
for which the quadrature rule is chosen to be Gaussian quadrature.

\subsection{ Convergence of 2D reconstruction algorithm} 

Let us consider the uniform norm $\|\cdot\|_\infty$ of continuous functions
on $B^2$. Convergence in the uniform norm guarantees the pointwise 
convergence of the reconstruction. First we give a formula for the operator 
norm.  Let us  denote by
$$
\psi_j:= \theta_{j,2m} = j\pi/(2m+1), \qquad 0 \le j \le 2m.
$$

\begin{prop} \label{prop:4.1}
Let $\|\CA_{2m}\|_\infty$ denote the operator norm of $\CA_{2m}$ as an 
operator from $C(B^2)$ into $C(B^2)$. Then
\begin{equation} \label{eq:Lambda}
 \|\CA_{2m}\|_\infty = \max_{(x,y)\in B^2} \Lambda_m(x,y), \qquad 
    \Lambda_m(x,y) := \sum_{\nu = 0}^{2m}\sum_{j=1}^{2m}\sin \psi_j
        |T_{j,\nu}(x,y)|. 
\end{equation} 
\end{prop}

\begin{proof}
By the definion of $\CR_{\phi} (f;t)$ we evidently have
$$
    |\CR_{\phi} (f;t)| \le \sqrt{1-t^2} \|f\|_\infty, 
         \quad 0 \le \phi \le 2\pi, \quad -1 < t < 1,
$$
as seen from the second equality of  \eqref{eq:Radon}. 
It follows immediately that 
\begin{equation}\label{eq:upperbd}
  |\CA_{2m}(f;x,y)| \le \|f\|_\infty \Lambda_m (x,y), 
     \quad\hbox{where}\quad (x,y) \in B^2.  
\end{equation}
Taking maximum in both side proves that $\|\CA_{2m}\|_\infty \le 
\max_{(x,y)\in B^2} |\Lambda_m(x,y)|$. To show that the equality holds,
let $(x_0,y_0)$ be a point in $B^2$ at which $\Lambda_m(x,y)$ attains 
its maximum over $B^2$. Recall that $I(\theta,t)$ denote a line segment
\eqref{eq:line} inside $B^2$. Let $\Sigma$ denote the set of intersection 
points of any two line segments $I(\phi_\nu,\cos \phi_j)$ and
$I(\phi_\mu,\cos \psi_i)$, 
$$
\Sigma : = \{(x,y): I(\phi_{\nu}, \cos \psi_j) \cap 
   I(\phi_{\mu}, \cos \psi_i), \quad (i,\mu) \ne (j,\nu)\}. 
$$
The set contains only finitely many points. Let $\varepsilon > 0$ be small 
enough so that a disk centered at a point in $\Sigma$ of radius $\varepsilon$
contains no other points in $\Sigma$. Let $\Sigma_\varepsilon$ denote the 
union of all such $\varepsilon$ disks. We construct a function $f_\varepsilon
\in C(B^2)$ as follows: 
$$
 f(x,y) = \sin \psi_j \sign T_{j,\nu}(x_0,y_0), \qquad (x,y) \in 
   I(\phi_\nu,\cos \phi_j)\setminus (I(\phi_\nu,\cos \phi_j)\cap\Sigma)
$$ 
for all $j,\nu$ and $\|f_\varepsilon\|_\infty =1$. Then 
$\CR_{\phi_\nu}(f_\varepsilon, \cos\psi_j) = \sign T_{j,\nu}(x_0,y_0) +
c_{j,\nu} \varepsilon$ for some costant $c_{j,\nu}$. Since there are only 
finitely many points in $\Sigma$, this shows that 
$$
 \|\CA_{2m}\|_\infty \ge |\CA_{2m} f_\varepsilon (x_0,y_0)|
 = \Lambda_m(x_0,y_0) - c \,\varepsilon 
 = \max_{(x,y)\in B^2} \Lambda_m(x,y) - c \,\varepsilon,  
$$
where $c$ is a constant that depends on $m$. Taking $\varepsilon \to 0$ 
completes the proof. 
\end{proof}

In the following we shall use the notation $A\approx B$  to mean that there 
exist two positive constant $c_1$ and $c_2$ such that $c_1 A \le B \le c_2 A$. 

Recall that $\CA_{2m} f$ is obtained from the partial sum $S_{2m} f$ of the
Fourier orthogonal expansion, upon using Gaussian quadrature. It is proved 
in \cite{X01} that the operator norm for the partial sum operator $S_{2m}$ is 
\begin{equation}
  \|S_{2m}\|_\infty  \approx m. 
\end{equation}
One would expect that the operator norm of $\CA_{2m}$ is worse than 
$\CO(m)$
due to the additional step of using Gaussian quadrature.  Our main result in
this section shows that the norm of $\CA_{2m}$ does not grow much worse.  

\begin{thm} \label{thm:A-norm}
For $\CA_{2m}$ defined in \eqref{eq:AlgoG},   
$$
     \|\CA_{2m} \|_{\infty} \approx m \log (m+1). 
$$
\end{thm}

This shows that the price we pay for using Gaussian quadrature to get  
to $\CA_{2m}$ is just a  $\log (m+1)$ factor.  This theorem will be proved in
the following subsection. 

For $f \in C(B^2)$, the quantity $E_n(f)_\infty$ defined in \eqref{eq:bestapp}
denotes the error of best approximation of $f$ by polynomials of degree at 
most $n$. It is proved in \cite{X05} that if $f \in C^{2r} (B^2)$, $r \in\NN$,
then 
\begin{equation} \label{eq:Enf}
        E_n (f) \le c n^{-2 r} \|\CD^r f \|_\infty,  \qquad n \ge 0,
\end{equation}
where $\CD$ is a second order partial differential operator defined in 
\eqref{eq:CD}.

As a consequence of Theorem \ref{thm:A-norm}, we have the following
result: 

\begin{thm} 
If $f \in C^2(B^2)$, then $\CA_{2m} f$ converges to $f$ uniformly. In fact, 
let $r$ be a positive integer; then for $f \in C^{2r}(B^2)$, 
$$ 
\|\CA_m f - f\|_\infty\le c \frac{\log (m+1)}{m^{2r-1}} \|\CD^r f \|_\infty. 
$$ 
\end{thm} 

\begin{proof}
Let $p$ be the best approximation for $f$ from $\Pi_{2m}^2$. By the 
definition of the operator norm, the fact that $\CA_{2m}$ preserves 
polynomials of degree up to $2m-1$, and the triangle inequality we see that 
\begin{align*}
 \|\CA_{2m} f - f \|_\infty  & \le \|CA_{2m}(f-p)\|_\infty + \|f-p\|_\infty \\
       &  \le (1+ \|\CA_{2m}\|_\infty ) E_{2m-1} (f)_\infty 
           \le c m \log (m+1) E_{2m-1} (f)_\infty 
\end{align*}
from which the stated inequality follows from \eqref{eq:Enf}.  
\end{proof}

This theorem shows that the Algorithm  \ref{algo:2D-G} does converge 
whenever $f $ is a $C^2(B^2)$ function. In other words, if the original image
is $C^2$ smooth then the reconstruction algorithm \ref{algo:2D-G} 
converges to the image pointwisely and uniformly. The speed of the
convergence depends on the smoothness of the function. 
 
The algorithm with multiplier function will likely have better convergence
behavior (recall Proposition \ref{prop:S-eta}), but the estimate is more 
difficult to establish. We will  report results along this line in future 
communications.
 
\subsection{ Proof of Theorem \ref{thm:A-norm}, lower bound} 

First we need a compact formula for the functions $T_{j,\nu}$ in 
\eqref{eq:TjvG}.  The notation
$$
  \cos  \theta_\nu(x,y) =  x \cos  \phi_\nu + y \sin \phi_\nu, \qquad 
        \phi_\nu = 2\pi \nu/(2m+1). 
$$
will be used throughout the rest of this paper. 

\begin{prop} 
For $0 \le \nu \le 2m$ and $1 \le j \le 2m$, 
\begin{align} \label{eq:Tjv}
(2m+1)^2T_{j,\nu}(x,y) = & \frac{- \sin \psi_j 
   [ 1- (-1)^j T_{2m+1}(\cos \theta_\nu(x,y)) ]}
      {2 ( \cos \theta_\nu(x,y) - \cos \psi_j)^2} \\ 
&  - (-1)^j \sin \psi_j \frac{(2m+1)U_{2m}( \cos \theta_\nu(x,y))} 
       {2  (\cos \theta_\nu(x,y) - \cos \psi_j)}. \notag
\end{align}
\end{prop}

\begin{proof}
To derive the formula, we start with the elementary trigonometric identity: 
$$
  \sum_{k=0}^{2m} (k+1) \cos((k+1) \theta) 
     = \frac{-1+(2m+2)\cos((2m+1)\theta) - (2m+1)\cos((2m+2)\theta)}
         {4 \sin^2(\theta/2)}.
$$
Let $\theta = \theta_\nu(x,y)$ in this proof.  We apply the above identity to 
\begin{align*}
&\sin \theta (2m+1)^2 T_{j,\nu}(x,y) = 
 \sum_{k=0}^{2m} (k+1) \sin((k+1)\psi_j ) \sin((k+1)\theta)\\
& \qquad\qquad  = \frac{1}{2}  \sum_{k=0}^{2m} (k+1) 
   \left[\cos((k+1)(\theta-\psi_j) )- \cos((k+1)(\theta+\psi_j))\right]
\end{align*}
and combine the two terms together as one fraction. The denominator of
the fraction is 
$$
2\cdot 4 \sin^2 \frac{\psi_j+\theta}{2}\sin^2 \frac{\psi_j-\theta}{2} 
    = 2(\cos \theta - \cos\psi_j)^2
$$
and the numerator of the fraction is 
\begin{align*}
 N(\theta):= \sin^2 \frac{\theta+\psi_j}{2} h_{j,m}(\theta-\psi_j)
    - \sin^2 \frac{\theta-\psi_j}{2} h_{j,m}(\theta+\psi_j)
\end{align*}
where 
$$
  h_{j,m}(\theta) = -1+(2m+2)\cos((2m+1)\theta) - (2m+1)\cos((2m+2)\theta).
$$
We write $N(\theta)$ as a sum of two terms, 
$$
N(\theta) =  N_1(\theta) - (2m+1) N_2(\theta).
$$ 
For the first term we use the identities 
\begin{align*}
 \cos ((2m+1) (\theta\pm \psi_j))  = (-1)^j \cos ((2m+1) \theta), 
\end{align*}  
and $2 \sin^2 (\theta/2) = 1- \cos \theta$ to get
\begin{align*}
  N_1(\theta):=&  \sin^2 \frac{\theta+\psi_j}{2} 
      ( -1+(2m+2)\cos [(2m+1)(\theta-\psi_j)] ) \\
    &  -    \sin^2 \frac{\theta-\psi_j}{2}  
      ( -1+(2m+2)\cos [(2m+1)(\theta+\psi_j)] ) \\
    = & ( -1+(2m+2)(-1)^j \cos [(2m+1)\theta] )
         (\cos(\theta-\psi_j) - \cos(\theta+\psi_j) )/2 \\      
    = & ( -1+(2m+2)(-1)^j \cos [(2m+1)\theta] )
         \sin \theta \sin\psi_j). 
\end{align*}
For the second term, we use 
$$
 \cos ((2m+2) (\theta \pm \psi_j))  = (-1)^j 
     \cos((2m+2)\theta \pm \psi_j),
$$
$2 \sin^2 (\theta/2) = 1- \cos \theta$ and the addition formula of the 
cosine function to obtain 
\begin{align*}
& N_2(\theta):=  \sin^2 \frac{\theta+\psi_j}{2} 
      \cos ((2m+2)(\theta-\psi_j)) - \sin^2 \frac{\theta-\psi_j}{2} 
      \cos((2m+2)(\theta+\psi_j)) \\
&  = \frac{ (-1)^j}{2} [
         \cos((2m+2)\theta-\psi_j) - \cos((2m+2)\theta+\psi_j)] \\
&  + \frac{(-1)^j}{2} [\cos (\theta-\psi_j) \cos((2m+2)\theta-\psi_j)
       -  \cos (\theta+\psi_j) \cos((2m+2)\theta+\psi_j)] \\   
&  =  (-1)^j \sin \psi_j \sin((2m+1)\theta) 
  + \frac{(-1)^j}{2} [-\cos (\theta-\psi_j) \sin (\theta+\psi_i)
        \sin ((2m+1)\theta) \\
   & \qquad \qquad \qquad \qquad \qquad \qquad \qquad \qquad
 + \cos (\theta+\psi_j) \sin (\theta-\psi_i)
        \sin((2m+1)\theta) ] \\
&  =  (-1)^j \sin \psi_j \sin ((2m+2)\theta) 
      + (-1)^{j+1} \sin ((2m+1)\theta)  \sin \psi_j \cos \psi_j, 
\end{align*}
where we have used the double angle formula for sine in the last step. 
Using the addition formula 
$\sin((2m+2)\theta) = \sin((2m+1)\theta) \cos \theta + 
\cos((2m+1)\theta) \sin \theta$, we obtain
$$
  N_2(\theta) = (-1)^j \sin \psi_j
     [\sin ((2m+1)\theta) (\cos \theta - \cos \psi_j) + 
       \cos ((2m+1)\theta) \sin \theta ].
$$
Putting the two terms together we obtain
\begin{align*}
 N(\theta)  = &  - \sin \theta \sin\psi_j + (-1)^j \sin \psi_j
        [(2m+2) \sin \theta \cos((2m+1)(\theta)) \\
      & - (2m+1)(\sin ((2m+1)\theta) (\cos \theta - \cos \psi_j) + 
       \cos((2m+1)\theta) \sin \theta )] \\
 = &  - \sin \theta \sin\psi_j [ 1- (-1)^j \cos ((2m+1)\theta) ] \\
      & - (2m+1) (-1)^j \sin \psi_j \sin ((2m+1)\theta) 
      (\cos \theta - \cos \psi_j). 
\end{align*}
Consequently, we have proved that 
\begin{align*}
 \sin \theta (2m+1)^2 T_{j,\nu}(x,y)  
 = & \frac{- \sin \theta \sin\psi_j [ 1- (-1)^j \cos ((2m+1)\theta) ]}
      {2 (\cos \theta - \cos \psi_j)^2} \\ 
&  - \frac{(2m+1) (-1)^j \sin \psi_j \sin ((2m+1)\theta)} 
       {2 (\cos \theta - \cos \psi_j)}.
\end{align*}
Hence, using the fact that $T_n(\cos \theta) = \cos n \theta$ and
$U_n(\cos \theta) = \sin(n+1)\theta / \sin \theta$, we conclude that
\begin{align*}
 & (2m+1)^2 T_{j,\nu}(x,y)  \\
 & =  \frac{- \sin \psi_j [ 1- (-1)^j T_{2m+1}(\cos \theta) ]}
      {2 ( \cos \theta - \cos \psi_j)^2}  
  - \frac{(2m+1) (-1)^j \sin \psi_j U_{2m}(\cos \theta)} 
       {2 (\cos \theta - \cos \psi_j)},
\end{align*}
which completes the proof.
\end{proof}

Throughout the rest of this section and in the following section, we will use 
the convention that $c$ denotes a generic constant, independent of $f$ and $m$,
its value may change from line to line. The elementary  facts 
$$   
\frac{2}{\pi} t \le \sin t\le t, \quad\hbox{for $0 \le t \le \frac{\pi}{2}$}, 
\quad \hbox{and}\quad 
    \cos \alpha - \cos\beta = 2 \sin \left(\tfrac{\beta - \alpha}{2}\right)
   \sin \left( \tfrac{\alpha+\beta}{2} \right)
$$
will be used repeatedly without further mention.

We now use the expression \eqref{eq:Tjv} to derive the lower bound for the
estimate.

\begin{prop}
$$
  \|\CA_{2m} \|_\infty \ge \Lambda_m \left(\cos \frac{\pi}{4m+2}, 
\sin \frac{\pi}{4m+2}\right) \ge c\,m\log(m+1).
$$
\end{prop} 

\begin{proof}
Let $x = \cos \frac{\pi}{4m+2}$ and $y =  \sin \frac{\pi}{4m+2}$. Then 
$\cos \theta_\nu(x,y) = \cos  \frac{(2\nu - 1/2)\pi}{2m+1}$, so that 
$\sin (2m+1) \theta_\nu(x,y) = 1$ and $\cos (2m+1) \theta_\nu(x,y) =0$. Let 
$\theta_\nu = \frac{(2\nu - 1/2)\pi}{2m+1}$. Then by \eqref{eq:Tjv}, 
$$
 (2m+1)^2  T_{j,\nu}(x,y) = - \frac{\sin \psi_j}{2} 
    \left[ \frac{1}{(\cos\theta_\nu - \cos \psi_j)^2} +
        \frac{2m+1}{\sin \theta_\nu (\cos\theta_\nu - \cos \psi_j)} \right]. 
$$
We will use the fact that $0 < \theta_\nu < \pi$ and $0 < \psi_j < \pi/2$
for $0 \le \nu, j \le m$.   Furthermore, for $0 < \theta_\nu \le \pi/2$ and 
$0 < \psi < \pi/2$ we have
$$
 \frac{\sin \psi_j}{ \sin \frac{\psi_j+\theta_\nu} {2} }
 = \frac{2 \sin \frac{\psi_j}{2}
  \cos \frac{\psi}{2} }{ \sin    \frac{\psi_j+\theta_\nu} {2}} 
  \le  2 \cos \frac{\psi_j}{2} \le 2,
$$
which also holds for $\pi/2 \le \theta_\nu \le \pi$ and $0 < \psi < \pi/2$, 
since then $\pi/4 \le \frac{\psi_j + \theta_\nu}{2} \le 3 \pi/4$ and 
$\sin \frac{\psi_j+\theta_\nu}{2} 
\ge \sqrt{2}/{2}$. Hence,  we have 
 \begin{align*}
 \frac{1} { (2m+1)^2 } \sum_{\nu = 0}^m \sum_{j=1}^m 
              \frac{ \sin^2 \psi_j } {(\cos\theta_\nu - \cos \psi_j    )^2 } 
           &  \le  \frac{1} { (2m+1)^2 } \sum_{\nu = 0}^m \sum_{j=1}^m 
              \frac{1} { \sin^2\frac{\theta_\nu - \psi_j}{2}  } \\
            & \le  \sum_{\nu = 0}^m \sum_{j=1}^m 
                \frac{1} {( 2 \nu  - j -1/2)^2}  \le c\, m.   
 \end{align*}
Consequently, it follows that 
\begin{align*}
  \Lambda_m(x,y) & \ge  \frac{1} { (2m+1)^2 } \sum_{\nu = 0}^m \sum_{j=1}^m 
   \sin \psi_j |T_{j,\nu}(x,y)| \\
    &  \ge  \frac{1} { 2m+1  } \sum_{\nu = 0}^m \sum_{j=1}^m 
 \frac{\sin^2 \psi_j} {|\sin \theta_\nu (\cos\theta_\nu-\cos\psi_j)|}-c\, m.
\end{align*}
The sum in the last expression is bounded below by
\begin{align*}
  &   \frac{8}{\pi^2} \frac{1} { 2m+1 } \sum_{\nu = 1}^m \sum_{j=1}^m 
          \frac{ \psi_j ^2} {  \theta_\nu |\theta_\nu - \psi_j|
  (\theta_\nu + \psi_j) } \\
  &  =   \frac{8}{\pi^3} \sum_{\nu = 1}^m \sum_{j=1}^m 
          \frac{ j ^2} { (2\nu - 1/2) |2\nu - j - 1/2| (2\nu +j -1/2)} \\
  & \ge   \frac{1}{\pi^3} \sum_{\nu = 1}^{m/2} \sum_{j=\nu}^{2\nu-1}
          \frac{1} {2\nu - j  }  = \frac{1}{\pi^3}  \sum_{\nu=1}^{m/2}
         \sum_{j= 1 }^{\nu} \frac{1} {j}   = \frac{1}{\pi^3}  \sum_{j=1}^{m/2}
               \frac{m/2 - j+1} {j} \\
  &  \ge \frac{1}{\pi^3}\left (\frac{m}{2} + 1 \right )
  \sum_{j=1}^{m/2} \frac{1}{j}
     - \frac{m}{2 \pi^3} \ge   c \, m \log (m+1). 
\end{align*}
This completes the proof. 
\end{proof} 

\subsection{Proof of Theorem \ref{thm:A-norm}, upper bound} 
We will also use the expression \eqref{eq:Tjv} to estimate 
$\Lambda_m(x,y)$ in \eqref{eq:Lambda} from above for all  
$(x,y) \in B^2$.  In the following we write $\Lambda(x,y) = \Lambda_m(x,y)$.

It is easy to see that the function $\Lambda(x,y)$ is invariant under the 
dihedral group $I_{2m+1}$; that is, it is invariant under the rotation of a
angle $\phi_\nu$ for $\nu = 0, 1, \ldots, 2m$. Hence, it suffices if we 
establish the estimate assuming $(x,y)$ is in the wedge 
$$
  \Gamma_m: = \{(x,y): x= r\cos \phi, y = r \sin \phi, \,\, r \ge 0,  \,   
     |\phi| \le  \varepsilon_m\} , \quad \varepsilon_m = \frac{\pi /2}{2m+1}. 
$$
Note that the set $\Gamma_m$ is symmetric with respect to the $y$ axis.

We start with a number of reductions.  The fact that $\phi_{2m+1-\nu} = 
2\pi - \phi_\nu$ shows  $\cos \phi_{2m+1-\nu} = \cos \phi_\nu$ and 
$\sin \phi_{2m+1-\nu} = - \sin \phi_\nu$, which implies that  
$\theta_{2m+1-\nu}(x,y) = \theta_\nu(x,-y)$.  Hence, using that 
$\theta_0(x,y) = x$, we can write
\begin{align*} 
  \Lambda (x,y) & =   \sum_{j=1}^{2m} \sin  \psi_j 
    \left[|T_{0,\nu}(x,y)| + \sum_{\nu=1}^{m}  
               \left( |T_{j,\nu} (x,y)| + |T_{j,\nu}(x,-y)| \right)  \right] \\
       & := H_{0} (x,y)+ \sum_{ \nu =1}^m \left( H_\nu (x,y)+  H_\nu (x,-y)
   \right). \notag       
\end{align*}
Since $(x,y) \in \Gamma_m$, we only need to consider the sum over 
$H_\nu(x,y)$.

The equation \eqref{eq:Tjv} shows that $T_{j,\nu}$ is naturally split
as a sum of two functions, 
\begin{align}
 &  T_{j,\nu}^{(1)} : = \frac{- \sin \psi_j 
   [ 1- (-1)^j T_{2m+1}(\cos \theta_\nu(x,y)) ]}
      {2 (\cos \theta_\nu(x,y) - \cos \psi_j)^2}  \label{eq:Tjv1}\\
  & T_{j,\nu}^{(2)} : =  - (-1)^j \sin \psi_j \frac{(2m+1) 
    U_{2m}( \cos \theta_\nu(x,y))} 
       {2 ( \cos \theta_\nu(x,y) - \cos \psi_j)}.\label{eq:Tjv2}
\end{align}
We shall denote the corresponding splitting of $H_{\nu}(x,y)$ by
$H_{\nu}^{(i)}(x,y)$ and the splitting of $\Lambda(x,y)$ (only the sum over
$H_\nu(x,y)$ by $\Lambda^{(i)}(x,y)$. Thus, 
$$
\Lambda^{(i)}(x,y):= \sum_{\nu=1}^m H_\nu^{(i)} (x,y)
    \quad\hbox{and}\quad 
H_\nu^{(i)} (x,y) :=   \sum_{j = 1}^{2m} \sin \psi_j|T_{j,\nu}^{(i)} (x,y)|, 
    \quad i = 1,2. 
$$

Next we split $H_\nu^{(i)}$ as two sums; the first one is over $j =1,2,
\ldots,m$ and the second one is over $j=m+1,m+2,\ldots, 2m$. Since 
$\psi_{2m+1-j} = \pi - \psi_j$, we have  $\cos \psi_{2m+1-j} = - 
\cos \psi_j$ so that we can write 
\begin{equation} \label{eq:Hv}
  H_\nu^{(i)} (x,y)  = H_{\nu,1}^{(i)}(x,y) + H_{\nu,2}^{(i)}(x,y),
\end{equation}
where, using that $1- (-1)^j \cos(2m+1)\theta = 
1- \cos (2m+1)(\theta -\psi_j)$,  
\begin{align} \label{eq:Hnu1}
\begin{split}
   H_{\nu,1}^{(1)} & = \frac{1}{(2m+1)^2} \sum_{j=1}^m \sin^2 \psi_j 
    \left | \frac{1- \cos (2m+1)(\theta -\psi_j)}
        {2(\cos \theta_\nu(x,y)-\cos \psi_j)^2}   \right|, \\
   H_{\nu,1}^{(1)} & = \frac{1}{(2m+1)^2} \sum_{j=1}^m \sin^2 \psi_j 
    \left | \frac{1- \cos (2m+1)(\theta -\psi_j)}
         {2(\cos \theta_\nu(x,y)+\cos \psi_j)^2}   \right|, 
\end{split}
\end{align}
and, using the fact that $(-1)^j \sin(2m+1) \theta =
 \sin (2m+1)(\theta -\psi_j)$,  
\begin{align} \label{eq:Hnu}
\begin{split}
H_{\nu,1}^{(2)}(x,y)  & =   \frac{1}{2m+1}  
 \sum_{j = 1}^{m} \sin^2 \psi_j  \left |
   \frac{\sin (2m+1)(\theta_\nu(x,y)-\psi_j)}
   {2 \sin \theta_\nu(x,y)(\cos \theta_\nu(x,y)-\cos \psi_j)} \right |, \\
H_{\nu,2}^{(2)}(x,y)  & =   \frac{1}{2m+1}  \sum_{j = 1}^{m} \sin^2 \psi_j   
 \left |  \frac{\sin (2m+1)(\theta_\nu(x,y)-\psi_j)} 
  {2\sin\theta_\nu(x,y) (\cos \theta_\nu(x,y) + \cos \psi_j)} \right |.
\end{split}
\end{align}
Let us denote the corresponding split of $\Lambda^{(i)}$ by $\Lambda_j^{(i)}$; 
that is, 
$$
\Lambda^{(i)}:= \Lambda^{(i)}_1 +  \Lambda^{(i)}_2, \quad\hbox{where}\quad
  \Lambda^{(i)}_j:= \sum_{\nu=1}^m  H_{\nu,j}^{(i)}, \quad i,j=1,2.  
$$
We only need to estimate one of the two terms. To see this, let us define
$$
   \wt \phi_\nu := (2\nu - 1) \pi /(2m+1), \qquad  1 \le \nu \le m. 
$$
Then $\cos \phi_{m+1-\nu} = \cos (\pi - \wt \phi_\nu) = - \cos \wt \phi_\nu$ 
and $\sin \phi_{m-\nu+1}  = \sin \wt \phi_\nu$. Consequently, 
$$
   \cos \theta_{m-\nu+1}(x,y) = - x \cos \wt \phi_\nu + y \sin \wt \phi_\nu = 
     - \cos \wt \theta_\nu (x, -y),
$$
where $\wt \theta_\nu (x, y) = x \cos \wt \phi_\nu + y \sin \wt \phi_\nu$. We 
will also use the notation $\wt H_{\nu,1}^{(i)}$ and $\wt \Lambda^{(i)}_j$
when $\theta_\nu(x,y)$ is replaced by $\wt \theta_\nu(x,y)$ in 
$H_{\nu,1}^{(i)}$.  It then follows from \eqref{eq:Hnu} that 
$$
H_{m-\nu+1,2}^{(i)}(x,y) =  \wt H_{\nu,1}^{(i)}(x, -y),  \qquad 1 \le \nu 
\le m, 
$$ 
and, consequently, $\Lambda^{(i)}_2(x,y) = \wt \Lambda^{(i)}_1(x,-y)$. Hence, 
the estimate for $\Lambda_2^{(i)}$ will be similar to the estimate for
$\Lambda_1^{(2)}$. In fact, set 
$$
\wt \Gamma_m: = \{(x,y): x= r\cos \phi, y = r \sin \phi, \,\, r \ge 0,  \,   
|\pi - \phi| \le  \varepsilon_m\} , \quad \varepsilon_m = \frac{\pi /2}{2m+1}; 
$$
then the estimate of $\Lambda^{(i)}_2(x,y)$ over $\wt \Gamma_m$ will be 
exactly the same as the estimate of $\Lambda^{(i)}_1(x,y)$ over $\Gamma_m$.
Thus, we only need to estimate one sum, which we choose to be 
$\Lambda_1^{(i)}$. 

We use the equations \eqref{eq:Hnu1} and \eqref{eq:Hnu} to carry out the 
estimate. The inequality 
\begin{equation} \label{eq:Un}
   |U_n(\cos t)| = \left |\frac{\sin(n+1) t}{\sin t} \right | \le n+1, 
   \qquad 0 \le t \le \pi,
\end{equation}
will be used several times.  The estimate is divided into several cases. 

\begin{lem} There exist constants $c_1$ and $c_2$ such that 
$$
 H_0^{(1)} (x,y) \le c_1 \quad\hbox{and}\quad 
 H_0^{(2)} (x,y)  \le c_2 \left(m +\frac{1}{2} \right)\log (m+1)
$$
for $(x,y) \in \Gamma_m$, where $c_1 < \pi^2 (2+\pi^2/12)$ and 
$c_2 <\pi^2/2+1$.  
\end{lem}

\begin{proof}
Since  $\cos \theta_0(x,y) = x$ and $(x,y) \in \Gamma_m$, we can write
$\theta =  \theta_\nu(x,y)$ with $0 \le \theta \le \pi/2$.  

We estimate the term $H_\nu^{(2)}$ first.  Since $\cos \psi_j > 0$ for 
$1 \le j \le m$ and $\cos \theta \ge 0$, it follows from \eqref{eq:Un} that 
\begin{align*}
 H_{0,2}^{(2)}(x,y) \le & \sum_{j=1}^m \frac{\sin^2 \psi_j}{\cos \psi_j} 
     \le  \sum_{j=1}^m \frac{1}{\cos \psi_{m-j}} 
  =  \sum_{j=1}^m \frac{1}{\sin\frac{(j+1/2)\pi}{2m+1}}\\
 \le & \frac{\pi}{2} \sum_{j=1}^m \frac{2m+1}{(j+1/2)\pi} \le
   \left(m+\frac12\right) \log (m+1).
\end{align*}

For the term  $H_{0,1}^{(2)}$, we need to consider the position of $\theta$ 
relative to that of $\psi_j$, which is divided into several further cases. 

\medskip\noindent
(A) If $0\le \theta \le \varepsilon_m = \pi/ (4m+2)$,  then by \eqref{eq:Un}
\begin{align*}
 H_{0,1}^{(2)}(x,y) \le & \sum_{j=1}^m \frac{\sin^2 \psi_j}
     {4 \sin \frac{\psi_j - \theta}{2} \sin \frac{\psi_j + \theta}{2} }
    \le  \pi^2 \sum_{j=1}^m \frac{\psi_j^2}{\psi_j^2-\theta^2} 
  \le \pi^2 \sum_{j=1}^m \frac{\psi_j^2}{\psi_j^2-\varepsilon_m^2} \\  
   \le &   \pi^2  \sum_{j=1}^m \frac{j^2}{(j^2-1/4)} \le 
      \pi^2 \left(m+ \frac{1}{4}\sum_{j=1}^m \frac{1}{(j^2-1/4)}\right) 
     \le  \pi^2 \left(m+ \frac12 \right) .
 \end{align*}

\noindent
(B) If $\psi_l -  \varepsilon_m \le \theta \le \psi_l + \varepsilon_m $, where 
$1 \le l \le m$,  then using  the fact that 
$$
\left | \frac{\sin(2m+1) (\theta-\psi) } {\sin \frac{\theta - \psi_l}{2} } 
\right| \le 
 2  \left| \frac{\sin(2m+1) \frac{\theta-\psi}{2} } 
{\sin \frac{\theta - \psi_l}{2} } \right | 
  \le 2m+1
$$   
we obtain that 
\begin{align*}
 H_{0,1}^{(2)}(x,y) \le & \frac{\sin^2 \psi_l}{4\sin\theta \sin \frac{\theta +
 \psi_l}{2} }
    +  \left( \sum_{j=1}^{l-1} + \sum_{j=l+1}^m \right) \frac{\sin^2 \psi_j}
          {4 \sin \theta \sin \frac{|\psi_j - \theta|}{2} \sin \frac{\psi_j + 
\theta}{2} } \\
     \le & \frac{\pi^2}{8}  \frac{\psi_l^2}{\theta (\psi_l + \theta)}  +  
\frac{\pi^3}{8(2m+1)} 
      \left( \sum_{j=1}^{l-1} +  \sum_{j=l+1}^m \right) \frac{\psi_j^2}{ 
\theta |\theta - \psi_j|
           (\theta+\psi_j)}\\
    \le & \frac{\pi^2}{8}  \frac{l}{l-1/2}  +  \frac{\pi^3}{8(2m+1)\theta} 
        \left( \sum_{j=1}^{l-1} \frac{j}{l -j -1/2}  +  \sum_{j=l+1}^m  
\frac{j}{j-l-1/2} \right).
\end{align*}
The first sum in the last expression is bounded by 
$$
  \frac{\pi^3}{8(2m+1)\theta}    \sum_{j=1}^{l-1} \frac{j}{l -j -1/2} \le  
  \frac{\pi^2(l-1)}{8(l-1/2)}  \sum_{j=1}^{l-1} \frac{1}{l -j-1/2} 
   \le  \frac{\pi^2}{4}(\log m +2),
$$
and the second sum is bounded by  
\begin{align*}
 & \frac{\pi^3}{8(2m+1)\theta}    \sum_{j=l+1}^m  \frac{j}{j-l-1/2} 
     \le  \frac{\pi^2}{8(l-1/2)} \sum_{j=1}^{m -l}  \frac{j+l}{j-1/2}  \\ 
 &  \qquad \qquad \le  \frac{\pi^2 m }{4} \sum_{j=1}^{m -l} \frac{1}{j-1/2}  
    \le \frac{\pi^2 m }{4} (\log m + 2).
\end{align*}
Putting these estimates together and use \eqref{eq:Hv}, we complete the proof 
for $H_0^{(2)}$. 

The proof for $H_0^{(1)}$ does not need require the splitting argument. 
By definition and \eqref{eq:Tjv1}, 
\begin{align*}
   H_0^{(1)} (x,y) = & \frac{1}{(2m+1)^2} \sum_{j=1}^{2m} \sin^2\psi_j 
     \frac{1-\cos ((2m+1)(\theta - \psi_j))}
      {8 \sin^2 \frac{\theta -\psi_j}{2} \sin^2 \frac{\theta +\psi_j}{2} }\\
     = & \frac{1}{(2m+1)^2} \sum_{j=1}^{2m} \sin^2\psi_j 
      \frac{\sin^2\left((2m+1)\frac{\theta  - \psi_j}{2}\right)} 
         {4 \sin^2 \frac{\theta -\psi_j}{2} \sin^2 \frac{\theta +\psi_j}{2} }.
\end{align*}
Let $\theta$ be fixed and $|\theta - \psi_k| \le \varepsilon_m$ for some $k$. 
Then 
by \eqref{eq:Un},
\begin{align*}
H_0^{(1)} (x,y) & \le  \frac{\sin^2 \psi_k}{4 \sin^2 \frac{\theta +\psi_k}{2}}
 + \frac{1}{(2m+1)^2} \sum_{j \ne k}  
 \frac{\sin^2 \psi_j}{4 \sin^2 \frac{\theta -\psi_j}{2}  
            \sin^2  \frac{\theta +\psi_j}{2} } \\
 \le & \, \frac{\pi^2}{4} \left (1 +  \sum_{j=1}^{k-1} \frac{1}{(k -j -1/2)^2}
   + \sum_{j=k+1}^{m} \frac{1}{(j - k -1/2)^2}\right) \\
  \le &\,  \frac{\pi^2}{2} \left(4+ \frac{\pi^2}{6} \right). 
\end{align*}
This completes the proof. 
\end{proof}

Putting all pieces together, it is readily seen that the proof of the 
Theorem \ref{thm:A-norm} follows from the conclusion of the following 
lemma. 

\begin{lem} There exist constants $c_1$ and $c2$ such that 
$$
   \Lambda^{(1)}_1(x,y) \le c_1\, (m+1)
\quad\hbox{and}\quad   
   \Lambda^{(2)}_1(x,y) \le c_2\, (m+1) \log (m+1)
$$
for $(x,y) \in \Gamma_m$, where $c_1 < \pi^2/96 + 7\pi^2 /24 +1/2$ and 
$c_2 < (3/2) \pi^2+\pi+1$.
\end{lem}

\begin{proof}
Again we consider the estimate for the sum $\Lambda_1^{(2)}$ first, which is
more difficult than the estimate for $\Lambda_1^{(1)}$. Using the polar 
coordinates $x = r \cos \phi$ and $y = r \sin \phi$, we can write
$\cos \theta_\nu(x,y)  = r \cos (\phi -\phi_\nu)$. For $(x,y) \in \Gamma_m$,  
$|\phi| \le \varepsilon_m$, so that 
$$
 \frac {(2\nu -1/2)\pi}{2m+1} \le \phi_\nu - \phi  \le
     \frac {(2\nu +1/2)\pi}{2m+1},
$$
from which we conclude that  $\cos \theta_\nu(x,y) \ge 0$ for $1\le\nu\le m/2$
and $\cos \theta_\nu(x,y) \le 0$ for $m/2+1\le \nu \le m$. 
  
We shall break the sum in $\Lambda_1^{(i)}$ into two parts, 
$$
\Lambda_1^{(i)}(x,y) = L_1^{(i)}(x,y) + L_2^{(i)}(x,y), 
$$
where the first one has the sum taken over $1 \le \nu \le m/2$ and the 
second one has the sum taken over  $m/2+1 \le \nu \le m$.

 \medskip  \noindent
{\bf Case 1.}   We consider the estimate of $L_2^{(i)}(x,y)$.  Note that 
$\cos \theta_\nu(x,y) \le 0$ and $\cos \psi_j \ge 0$ for the terms in 
$L_2^{(i)}(x,y)$. If $\sin \theta_\nu (x,y) \ge  \sqrt{2}/2$,  then 
 \begin{align*}
 L_2^{(2)}(x,y) & = \frac{1}{2m+1} \sum_{\nu = m/2+1}^m \sum_{j=1}^m 
 \sin^2 \psi_j \left| \frac{ \sin((2m+1)\theta_\nu(x,y)) }{2\sin\theta_\nu(x,y)
             ( \cos \theta_\nu(x,y) - \cos \psi_j ) } \right| \\
             & \le \frac{1}{2m+1} \sum_{\nu = m/2+1}^m \sum_{j=1}^m 
                  \frac{1}{\sqrt{2}} \frac{\sin^2 \psi_j}{ \cos \psi_j  }   
              \le \frac{1}{4 \sqrt{2} }   \sum_{j=1}^m 
                  \frac{1}{ \sin \frac{(j+1/2)\pi }{2m+1}}  \\
                 & \le \frac{2m+1}{8 \sqrt{2} }  \sum_{j=1}^m \frac{1}{j+1/2} 
    \le \frac{1}{4\sqrt{2}} \left(m+\frac{1}{2}\right) \log (m+1). 
 \end{align*}
If $\sin \theta_\nu (x,y) \le  \sqrt{2}/2$,  then $-\cos \theta_\nu(x,y) =
 \cos (\pi -\theta_\nu(x,y)) \ge \sqrt{2}/2$. Furthermore, since 
$-\cos \theta_\nu(x,y) = - r \cos (\phi_\nu -\phi) \le  - \cos(\phi_\nu-\phi)$ 
for $m/2+1 \le \nu \le m$, we have $ \theta_\nu(x,y) \le \phi_\nu - \phi$  
and  $|\sin \theta_\nu(x,y)| = \sin (\pi - \theta_\nu(x,y)) \ge 
(2/\pi) (\pi - \phi _\nu + \phi)$. Hence, for $\sin \theta_\nu (x,y) \le 
\sqrt{2}/2$, 
\begin{align*}
  L_2^{(2)} (x,y) &  \le \frac{1}{2m+1} \sum_{\nu = m/2+1}^m \sum_{j=1}^m 
   \frac{1}{\sqrt{2}} \frac{2\sin^2 \psi_j}{|\sin \theta_\nu(x,y) |}  \\
 & \le \frac{\pi}{2 \sqrt{2} } \sum_{\nu = m/2+1}^m \frac{1}{\pi - 
       \phi_\nu+\phi} = 
      \frac{1}{2 \sqrt{2} } \sum_{\nu =1}^{m/2} \frac{2m+1}{2\nu +1} \\
 &    \le  \frac{1}{2 \sqrt{2} } \left(m+\frac{1}{2}\right)\log (m+1). 
\end{align*}
The case  $L_2^{(1)}(x,y)$ is again easier. We have 
 \begin{align*}
 L_2^{(1)}(x,y) & = \frac{1}{(2m+1)^2} \sum_{\nu=m/2+1}^m 
  \sum_{j=1}^m \sin^2\psi_j
    \left | \frac{\sin^2 \left((2m+1)\frac{\theta_\nu(x,y)-\phi_j}{2}\right) }
        {2 (\cos \theta_\nu(x,y) - \cos \psi_j)^2} \right| \\
   & \le \frac{m}{4 (2m+1)^2} \sum_{j=1}^m \frac{\sin^2 \psi_j} {\cos^2\psi_j} 
      \le  \frac{m}{16} \sum_{j=1}^m   \frac{1}{ (j+1/2)^2}  \le 
  \frac{\pi^2}{96} \, m
 \end{align*} 
  
\medskip\noindent  
{\bf Case 2.}  The estimate of $L_1^{(i)}(x,y)$ under the assumption that
$$
0 \le r \le \cos (\psi_m - \varepsilon_m) = \sin \pi/(2m+1),
$$
where $\cos \theta_\nu(x,y) = r \cos (\phi - \phi_\nu)$ in the polar 
coordinates for $(x,y)$.  The fact that $\cos \theta_\nu(x,y) \le r$ implies 
that $\theta_\nu(x,y) \ge 
\psi_m - \varepsilon_m = 
\pi/2 - 2\varepsilon_m$, so that $\sin \theta_\nu(x,y) > 1/2$. Furthermore, 
$\pi /4 \le \frac{\theta_\nu(x,y) + \psi_j}{2} \le 3\pi/4$, which implies that
$\sin \frac{\theta_\nu(x,y) +\psi_j} {2} \ge \sqrt{2}/2 >  1/2$. Consequently,
\begin{align*}
  L_1^{(2)}(x,y) & \le \frac{1}{2m+1} \sum_{\nu =1}^{m/2} \sum_{j=1}^m 
     \sin^2 \psi_j \left| \frac{ \sin ((2m+1) (\theta_\nu(x,y) -\psi_j)) }{
        2 \sin \frac{\theta_\nu(x,y) - \psi_j}{2} \sin \frac{\theta_\nu(x,y) 
  + \psi_j}{2}} \right| \\ 
      & \le \frac{1}{2m+1} \sum_{\nu =1}^{m/2} \sum_{j=1}^{m-1} 
         \frac{\sin^2 \psi_j} { \sin \frac{|\theta_\nu(x,y) - \psi_j|}{2}} + 
                       \sum_{\nu=1}^m \frac{\sin^2 \psi_m}
                           {\sin \frac{|\theta_\nu(x,y) +\psi_m|}{2}} \\
      & \le  \frac{ m \pi }{2(2m+1)} \sum_{j =1}^{m-1} 
             \frac{1}{ \pi/2 - \psi_j - 2 \varepsilon_m} + 2 m \\
      & \le  \frac{m}{2} \sum_{j=1}^{m-1} \frac{1}{ m -j - 1/2} + 2 m \le 
            \,m (\log m +2).             
 \end{align*}
Similarly, we have the estimate 
\begin{align*}
L_1^{(1)}(x,y) = &  \frac{1}{(2m+1)^2} \sum_{\nu =1}^{m/2} \sum_{j=1}^m 
         \sin^2 \psi_j \left| \frac{ \sin^2 \left( (2m+1) 
       \frac{\theta_\nu(x,y) -\psi_j}{2}\right) } {
        4 \sin^2\frac{\theta_\nu(x,y) - \psi_j}{2} 
           \sin^2 \frac{\theta_\nu(x,y) + \psi_j}{2}} \right| \\ 
           & \le \frac{\pi^2}{4(2m+1)^2} \sum_{\nu =1}^{m/2} \sum_{j=1}^{m-1} 
                   \frac{1} {(\theta_\nu(x,y) - \psi_j)^2 } + 
                 \sum_{\nu =1}^{m/2} \frac{ \sin^2 \psi_m }
                    {4 \sin^2 \frac{\theta_\nu (x,y) + \psi_m }{2}} \\
           & \le  \frac{m}{2}  \sum_{j=1}^{m-1} \frac{1}{ (m-j - 1/2)^2} 
              +  \frac{m}{2}  \le  \frac{m}{2}\left(\frac{\pi^2}{6} +5\right). 
\end{align*}

\noindent\medskip
{\bf Case 3.} The estimate of $L_1^{(i)}(x,y)$ under the assumption that
$$
r \ge \cos (\psi_m - \varepsilon_m) = \sin \pi/(2m+1),
$$ 
where 
$\cos \theta_\nu(x,y) = r \cos (\phi - \phi_\nu)$ in the polar coordinates for 
$(x,y)$.  Note that $\cos \theta_\nu(x,y) \ge 0$ for the terms in 
$L_1^{(i)}(x,y)$. In this case,  $\cos \theta_\nu(x,y) - \cos \psi_j$ can be 
zero. For a fixed $r$ let $l$ be the index such that 
 $$
 r = \cos \theta, \qquad  |\theta - \psi_l | \le \varepsilon_m, \qquad  
1 \le l \le m.
$$ 
We split  the summation over $j$ as two sums, one over $1 \le j \le l$ and 
the other over  $l+1 \le  j \le m$. This leads to a split of the sums in 
$L_1^{(i)}(x,y)$, 
$$
  L_1^{(i)}(x,y)  = \sum_{\nu = 1}^{m/2} \sum_{j=1}^l  +
  \sum_{\nu = 1}^{m/2} \sum_{j=l+1}^m  
      :=   L_{1,1}^{(i)}(x,y) + L_{1,2}^{(i)}(x,y).
$$
We  consider the two cases separately. 
 
\medskip\noindent
(A)   The estimate for $L_{1,1}^{(i)}$.  
For $j \le l$, $r = \cos\theta \le \cos (\psi_l - \varepsilon) < 
\cos \psi_j$, so that $\cos \theta_\nu(x,y) \le r < \cos \psi_j$, and 
consequently, 
$$
|\cos \theta_\nu(x,y)  - \cos \psi_j | \ge \cos \psi_j  - \cos \theta_\nu (x,y)
\ge  \cos \psi_j - \cos (\psi_l - \varepsilon_m) > 0.
$$
Furthermore, the fact that $\cos \theta_\nu(x,y) = r \cos (\phi-\phi_\nu)\le
  \cos (\psi_l-\varepsilon_m)$
 implies that $\theta_\nu(x,y) \ge  \psi_l - \varepsilon_m$. Therefore
 \begin{align*}
   L^{(2)}_{1,1} (x,y) \le &  \frac{1}{2m+1} \sum_{\nu=1}^{m/2} \sum_{j=1}^l   
       \frac{\sin^2 \psi_j  }{ 4 |\sin \theta_\nu(x,y)
 \sin \frac{\theta_\nu(x,y) - \psi_j}{2} 
         \sin   \frac{\theta_\nu(x,y) + \psi_j}{2}|} \\
   \le &   \frac{ \pi^3 m }{8(2m+1)} 
   \frac{1}{\psi_l - \varepsilon_m} \sum_{j=1}^l   
       \frac{\psi_j^2} {(\psi_l+\psi_j - \varepsilon_m)
(\psi_l-\psi_j -\varepsilon_m)} \\
   \le &   \frac{ \pi^2 m }{4(2l-1)} \sum_{j=1}^l   \frac{j}{ l-j-1/2}
   \le  \frac{ \pi^2 m }{2} \log m. 
\end{align*}
The estimate for $L_{1,1}^{(1)}$ is similar, 
\begin{align*}
   L^{(1)}_{1,1} (x,y) \le &  \frac{1}{(2m+1)^2} 
        \sum_{\nu=1}^{m/2} \sum_{j=1}^l   
       \frac{\sin^2 \psi_j  }{ 8 |\sin^2 \frac{\theta_\nu(x,y) - \psi_j}{2} 
          \sin^2  \frac{\theta_\nu(x,y) + \psi_j}{2}|} \\
    \le &  \frac{ \pi^2 m }{4} \sum_{j=1}^l  \frac{1}{( l-j-1/2)^2}
     \le \frac{ \pi^2 m }{16} \left( \frac{\pi^2}{6} + 4\right). 
\end{align*}

\medskip\noindent
(B)   The estimate for $L_{1,2}^{(i)}$.  For $l < j \le m$, $r=\cos\theta \ge 
\cos (\psi_l + \varepsilon) > \cos \psi_j$; hence,  
$\theta_\nu(x,y) - \cos \psi_j$ can be zero.  Since $\cos \theta_\nu (x,y) = 
r \cos (\phi_\nu - \phi)$ is bounded by both 
$r \le \cos (\psi_l - \varepsilon_m)$ and $\cos (\phi_\nu - \varepsilon_m)$
as $|\phi| \le \varepsilon_m$, it follows that 
$$
  \theta_\nu (x,y) \ge \max \{\psi_l, \phi_\nu\} - \varepsilon_m 
    : = z_{l,\nu} - \varepsilon_m > 0. 
$$
For each $\nu$, we choose an index $j_\nu$ such that 
$$
    |\theta_\nu (x,y) - \psi_{j_\nu} | \le \varepsilon_m. 
$$
We need an estimate on $j_\nu$. Since 
\begin{align*}
  2 \sin^2 \frac{\psi_{j_\nu} - \varepsilon_m}{2} = & 
    1 - \cos (\psi_{j_\nu} - \varepsilon_m)
     \le 1 - \cos (\theta_\nu(x,y)) = 1- r \cos (\phi_\nu - \phi) \\
    & = 1- r + r(1-\cos (\phi_\nu - \phi)) \le 2 \sin^2 \frac{\psi_l}{2} +
       2 \sin^2 \frac{\phi_\nu - \phi}{2}, 
\end{align*}\
we obtain  
$
 (\phi_{j_\nu} - \varepsilon_m)^2 \le  \frac{\pi^2}{4} 
    (\psi_l^2 + (\phi_\nu - \phi)^2), 
$
from which it follows that 
\begin{equation} \label{eq:jv}
 j_\nu - \frac{1}{2} \le \frac{\pi}{2} \sqrt{l^2 + (2\nu)^2} 
    \le \frac{\pi}{2} (l+ 2\nu).  
\end{equation} 

Using the fact that for $j \ne j_\nu$, 
$$
| \theta_\nu(x,y) - \psi_j | \ge  
 |\psi_{j_\nu} - \psi_j | - | \theta_\nu(x,y) - \psi_{j_\nu} |
    \ge |\psi_{j_\nu} - \psi_j| - \varepsilon_m, 
$$
we obtain the estimate 
\begin{align*}
   L_{1,2}^{(2)}(x,y) & =  \frac{1}{2m+1}  \sum_{\nu =1}^{m/2} \sum_{j=l}^m 
        \sin^2 \psi_j \left| \frac{ \sin ((2m+1)(\theta_\nu(x,y) -\psi_j)) }
   {4\sin \theta_\nu(x,y) 
         \sin \frac{\theta_\nu(x,y) - \psi_j}{2} \sin \frac{\theta_\nu(x,y)
   + \psi_j}{2}} \right | \\ 
  & \le  \sum_{\nu =1}^{m/2}  \frac{\sin^2 \psi_{j_\nu} } 
  {4 \sin \theta_\nu(x,y)   \sin \frac{\theta_\nu(x,y) + \psi_{j_\nu}}{2}}  \\
     & \qquad\quad 
          +  \frac{1}{2m+1}  \sum_{\nu =1}^{m/2} \sum_{j \ne j_\nu}
             \frac{ \sin^2 \psi_j }{4 \sin \theta_\nu(x,y) 
         \sin \frac{|\theta_\nu(x,y) - \psi_j|}{2} 
  \sin \frac{\theta_\nu(x,y) + \psi_j}{2}}\\ 
    & \le \frac{\pi^2}{8}  
       \sum_{\nu =1}^{m/2} \frac{\psi_{j_\nu}}{z_{l,\nu} - \varepsilon_m} 
         + \frac{\pi^2}{2(2m+1)} \sum_{\nu =1}^{m/2} 
             \frac{1}{z_{l,\nu} - \varepsilon_m} 
 \sum_{j \ne j_\nu} \frac{\psi_j}{ |\psi_{j_\nu} - \psi_j| - \varepsilon_m }. 
\end{align*}
Using the inequality \eqref{eq:jv}, the first sum in the last expression 
is bounded by 
$$
\frac{ \pi^3}{16} \sum_{\nu =1}^{m/2} \frac{ l + 2\nu} {\max\{l, 2\nu\} - 1/2}
  \le  \frac{ \pi^3}{8} m,
$$
while the second sum is bounded by, again using \eqref{eq:jv}, 
\begin{align*}
 &  \frac{\pi^2}{2(2m+1)} \sum_{\nu =1}^{m/2}
  \frac{1}{z_{l,\nu} - \varepsilon_m} 
            \left( \sum_{j = l}^{j_\nu-1}  \frac{j}{  j_\nu -j  - 1/2} + 
           \sum_{j = j_\nu+1}^{m}  \frac{j}{ j- j_\nu   - 1/2}   \right) \\
 &   \le   \frac{\pi}{2} \sum_{\nu =1}^{m/2} \frac{1}{ \max\{l, 2\nu\} - 1/2} 
            \left( j_\nu \log (j_\nu +1) + 
          \sum_{j = 1}^{m- j_\nu}  \frac{j + j_\nu}{ j - 1/2}   \right) \\
 &   \le \frac{ \pi}{2} \sum_{\nu =1}^{m/2} \frac{1}{ \max\{l, 2\nu\} - 1/2} 
      \left( j_\nu \log (j_\nu +1) + (m-j_\nu ) +( j_\nu +1/2)\log (m-j_\nu+1) 
                  \right) \\
 &   \le  \frac{\pi}{2} (\pi+1) 
      \log (m+1) \sum_{\nu =1}^{m/2} \frac{l + 2 \nu }{ \max\{l, 2\nu\} - 1/2} 
        +  m \sum_{\nu =1}^{m/2}  \frac{1 }{ \max\{l, 2\nu\} - 1/2}   \\
 &   \le (\pi^2+\pi+1))  m \log (m+1) . 
 \end{align*}
This completes the estimate for $L_{1,2}^{(2)}$. Similarly, but more easily, 
we have by \eqref{eq:jv},
\begin{align*}
 L_{1,2}^{(1)} (x,y)    & \le  \sum_{\nu =1}^{m/2}  \frac{\sin^2 \psi_{j_\nu}} 
             { 4\sin^2 \frac{\theta_\nu(x,y) + \psi_{j_\nu}}{2}}  
                +  \frac{1}{4 (2m+1)^2}  \sum_{\nu =1}^{m/2} \sum_{j \ne j_\nu}
                    \frac{1}{ \sin^2\frac{\theta_\nu(x,y) - \psi_j}{2}} \\ 
            &  \le  \frac{m}{2} +  \frac{m}{8} \sum_{j \ne j_\nu} 
           \frac{1}{ (j_\nu- j)^2 } \le 
          \frac{m}{2} +  \frac{m}{4} \sum_{j=1}^\infty \frac{1}{j^2}   
         =  \frac{m}{2} \left(1+ \frac{\pi^2}{12}\right).
\end{align*}
This competes the estimate of $L_{1,2}^{(i)}$. 

Putting all cases together completes the proof of the lemma. It is easy to 
see that the largest constants in front of $m \log (m+1)$ or $m$ in the three 
cases come from Case 3. 
\end{proof}
 
\medskip

\noindent
{\it Acknowledgment}: The author thanks Dr. Christoph Hoeschen and 
Dr. Oleg Tischenko for stimulating discussions about the computed tomography. 
He also thanks a referee for his careful reading and comments.


\begin{thebibliography}{99}

\bibitem{BX}
        B. Bojanov, Yuan Xu, 
        Reconstruction of a polynomial from its Radon projections, 
        \textit{SIAM J. Math. Anal.}, to appear. 

\bibitem{DX}
        C. F. Dunkl, Yuan Xu,
        \textit{Orthogonal polynomials of several variables},
        Cambridge Univ. Press, 2001. 
 
\bibitem{LS} 
      B. Logan, L. Shepp, 
      Optimal reconstruction of a function from its projections, 
      \textit{Duke Math. J.} \textbf{42} (1975), 645-659.    
       
\bibitem{Marr} 
        R. Marr, 
        On the reconstruction of a function on a circular domain from 
        a sampling of its line integrals, 
        \textit{J. Math. Anal. Appl.}, \textbf{45} (1974), 357-374.
 
\bibitem{N} 
        F. Natterer,
        \textit{The mathetmaics of computerized tomography},
        Reprint of the 1986 original. Classics in Applied Mathematics, 
        32. SIAM, Philadelphia, PA, 2001. 

\bibitem{NW} 
        F. Natterer, F. W\"ubbeling, 
        \textit{Mathematical Methods in Image Reconstruction},
        SIAM, Philadelphia, PA, 2001. 
 
 \bibitem{X99}
        Yuan Xu,
        Summability of Fourier orthogonal series for Jacobi weight on a 
        ball in $\RR^d$, \textit{Trans. Amer. Math. Soc.}  \textbf{351} (1999),
        2439-2458.

\bibitem{X00}
        Yuan Xu,
        Funk-Hecke formula for orthogonal polynomials on spheres and
        on balls, 
        \textit{Bull. London Math. Soc.} \textbf{32} (2000), 447-457.
 
\bibitem{X01}
        Yuan Xu,
        Representation of reproducing kernels and the Lebesgue constants
        on the ball,       
        \textit{J. Approx. Theory} \textbf{112} (2001), 295-310.

\bibitem{X05}
        Yuan Xu,
        Weighted approximation of functions on the unit sphere, 
        \textit{Const. Approx.} \textbf{21} (2005),  1-28. 

\end{thebibliography}
\end{document}